\documentclass{article}
\usepackage{amsfonts,amssymb,graphicx}

\renewcommand{\labelenumi}{(\arabic{enumi})} 
\newtheorem{lemma}{Lemma}[section]
\newtheorem{theorem}{Theorem}[section]
\newtheorem{corollary}{Corollary}[section]
\newcommand{\bproof}{\paragraph {\bf Proof. }}

\newcommand{\df}{\displaystyle\frac}
\newcommand{\ds}{\displaystyle}
\newcommand{\hsp}{\hspace{\parindent}}

\newcommand{\qed}{\square}

\newcommand{\RR}{{\mathbb{R}}}
\newcommand{\ZZ}{{\mathbb{Z}}}
\newcommand{\bw}{{\bf w}}

\newcommand{\be}{{\bf e}}
\newcommand{\bx}{{\bf x}}

\newcommand{\bc}{{\bf c}}

\newcommand{\bz}{{\bf z}}
\newcommand{\bv}{{\bf v}}
\newcommand{\bn}{{\bf n}}

\newcommand{\ra}{{\rightarrow}}
\newcommand{\sL}{{\mathcal{L}}}
\newcommand{\sC}{\mathcal{C}}
\newcommand{\sD}{{\mathcal{D}}}
\newcommand{\sG}{{\mathcal{G}}}
\newcommand{\sE}{{\mathcal{E}}}
\newcommand{\sH}{{\mathcal{H}}}

\newcommand\hl[1]{\mbox{ #1}}

\begin{document}
\title{ The number of Reidemeister Moves Needed for Unknotting}
\author{Joel Hass and Jeffrey C. Lagarias
\footnote{This paper grew out of work begun
while the authors were visiting the Mathematical
Sciences Research Institute in Berkeley in 1996/7. Research at MSRI is
supported in part by NSF grant DMS-9022140. The first author is
partially supported by NSF grant DMS-9704286.}}
 
\maketitle

\begin{abstract}
There is a positive constant $c_1$ such that for any diagram $\sD$ representing
the unknot, there is a sequence of at most $2^{c_1 n}$ Reidemeister moves that will
convert it to a trivial knot diagram, where $n$ is the number of crossings in $\sD$.
A similar result holds for elementary moves on a polygonal knot $K$ embedded in the 1-skeleton of the
interior of a compact, orientable, triangulated $PL$ 3-manifold $M$. There is a positive constant
$c_2$ such that for each $t \geq 1$, if $M$ consists of $t$ tetrahedra, and $K$ is unknotted,
then there is a sequence of at most $2^{c_2 t}$ elementary moves in
$M$ which transforms $K$ to a
triangle contained inside one tetrahedron of $M$.
We obtain explicit values for $c_1$ and $c_2$.
\end{abstract}

{\em Keywords}:
knot theory, knot diagram, Reidemeister move, normal surfaces, computational complexity
\newpage
\begin{center}
{\Large \bf The Number of Reidemeister Moves \\
\smallskip Needed for Unknotting} \\ \bigskip
{\large {\em Joel Hass}} \\ \smallskip
University of California, Davis \\
Davis, CA 95616 \\ \bigskip
{\large {\em Jeffrey C. Lagarias}} \\ \smallskip AT\&T Labs -- Research \\
Florham Park, NJ 07932 \bigskip \\

\vspace*{2\baselineskip}
(\today) \\
\vspace*{1\baselineskip}
\vspace{2\baselineskip}
\end{center}
\setlength{\baselineskip}{1.5\baselineskip}
\section{Introduction}
\hsp
A knot is an embedding of a circle $S^1$ in a 3-manifold $M$, usually taken to be
$\RR^3$ or $S^3$.
In the 1920's Alexander and Briggs \cite[\S4]{AleBri26} and Reidemeister
\cite{R1a} observed that questions about ambient isotopy of polygonal knots in
$\RR^3$ can be reduced to combinatorial
questions about knot diagrams. These are labeled planar graphs with overcrossings and
undercrossings marked, representing a projection of the knot onto a plane.
They showed that any ambient isotopy of a polygonal knot can be achieved by a finite sequence
of piecewise linear
moves which slide the knot across a single triangle, which are called 
{\em elementary moves} (or $\Delta$-moves). They also showed that two knots were ambient isotopic if and only if
their knot diagrams were equivalent under a finite sequence of local combinatorial changes, now
called {\em Reidemeister moves}, see \S7. 

A triangle in $M$ defines a {\em trivial knot}, and a loop in the plane with no crossings is
said to be a {\em trivial knot diagram}. A knot diagram $\sD$ is {\em unknotted} if it is
equivalent to a trivial knot diagram under Reidemeister moves.

We measure the complexity of a knot diagram $\sD$ by using its {\em crossing number}, the number of
vertices in the planar graph $\sD$, see \S7. A problem of long standing is to determine an upper
bound for the number of Reidemeister moves needed to transform an unknotted knot diagram $\sD$ to
the trivial knot diagram, as an explicit function of the crossing number $n$, see Welsh \cite[p.~95]{Wel93}. This paper obtains such a bound. 

\begin{theorem} \label{thm101}
There is a positive constant $c_1$, such that for each $n \geq 1$, any unknotted knot diagram
$\sD$ with $n$ crossings can be transformed to the trivial knot diagram using at most $2^{c_1 n}$
Reidemeister moves.
\end{theorem}

We obtain the explicit value of $10^{11}$ for $c_1$;
this value can clearly be improved.
In this paper we have not attempted to get an optimal value for
$c_1$, and several of our estimates are accordingly rather 
rough.

In 1934 Goeritz
\cite{Goe34} showed that there exist diagrams $\sD$ of the unknot such that any sequence of
Reidemeister moves converting $\sD$ to the trivial knot must pass through some intermediate knot
diagram $\sD'$ that has more crossings than $\sD$. 

\begin{figure}[hbtp]
\centering
\includegraphics[width=.5\textwidth]{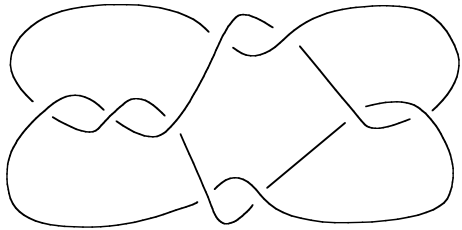}
\caption{ Goeritz's unknot.}
\label{goeritz}
\end{figure} 

It is not known how much larger the number of crossings needs to be for a
general $n$-crossing knot, but Theorem~\ref{thm101} implies the following upper
bound. 

\begin{corollary} \label{crossings.bound} 
Any unknotted knot diagram $\sD$ with $n$ crossings can be transformed by Reidemeister
moves to a trivial knot diagram through a sequence of knot diagrams each of which has at most $2^{c_1 n}$ crossings. \end{corollary}

This holds because each Reidemeister move creates or destroys at most two crossings, and all crossings must be eliminated by
the final Reidemeister move. It seems possible
that for each knot diagram representing the unknot there exists some
sequence of Reidemeister moves taking it to the trivial knot diagram which keeps the crossing
number of all intermediate knot diagrams bounded by a polynomial in $n$. 

In a general triangulated 3-manifold, a trivial knot is a triangle lying inside one
tetrahedron and an elementary move takes place entirely within a single tetrahedron of the
triangulation. The proof of Theorem~\ref{thm101}
is based on a general result that bounds the number of
elementary moves needed to transform a closed polygonal curve $K$ embedded in the 1-skeleton of
a triangulated $PL$ 3-manifold $M$ to a trivial knot.

\begin{theorem} \label{thm102}
There is a constant $c_2$ such that for each $t \geq 1$, any compact,
orientable, $PL$ 3-manifold with boundary ($M,
\partial M$) triangulated by
$t$ tetrahedra has the following property: If $K$ is an unknot embedded in the 1-skeleton of interior 
($ M$), then $K$ can be isotoped to the trivial knot using at most $2^{c_2 t}$ elementary moves.
\end{theorem}

We obtain the explicit value $c_2 = 10^7$. The hypothesis that $K$ lies in 
the interior of $M$ can be weakened to
requiring only that $K$ lies in $M$, as we indicate in \S6.

The methods we use are based on the normal surface theory
of Kneser \cite{Kne29} and Haken \cite{Hak61}.
Haken applied this theory to
obtain an algorithm to decide if a knot is trivial.
See \cite{Hass97} for a survey of algorithms to recognize unknotting, and
\cite{BM92} and \cite{BH97}
for possible alternate approaches based on braids.


The paper is organized as follows. In \S2 we define elementary moves and outline the proof of
Theorem~\ref{thm102}. In
\S3, \S4 and \S5 we bound the number of elementary moves needed to isotop a curve in three
distinct settings. First we consider an isotopy of a curve to a point across a compressing
disk in a 3-manifold, secondly we consider an isotopy through a solid torus from a core to a
longitude, and finally we consider an isotopy across a surface between two isotopic curves. In
\S6 we tie these results together to complete the proof of Theorem~\ref{thm102}. We then consider the
more special case of knots in $\RR^3$. In \S7 we outline the proof of Theorem~\ref{thm101} and
begin it by relating elementary moves to Reidemeister moves. In \S8 we start with a knot
diagram and show how to construct a triangulation of a convex submanifold of $\RR^3$ containing
the corresponding knot on its 1-skeleton. This allows us to apply Theorem~\ref{thm102}. In \S9 we
complete the proof of Theorem~\ref{thm101}. 

{\bf Note:} After completing this paper we learned that
S. Galatolo has obtained results about the case of Reidemeister moves in
$\RR^3$ similar to
those stated in Theorem~\ref{thm101} and Corollary~\ref{crossings.bound}, and at
about the same time.  Galatolo's constructions are also
based on an analysis of normal surfaces.
An announcement of his results will appear in \cite{Galatolo}.

\section{Elementary Moves in a 3-manifold} \hsp
In this section we define elementary moves for knots and links. Elementary moves are
well-defined in arbitrary triangulated $PL$ 3-manifolds. In section~7 we will discuss
Reidemeister moves, which are associated to knots and links in
$\RR^3$, and relate these two concepts for manifolds $PL$-embedded in $\RR^3$. 

A (piecewise linear, unoriented) {\em knot} $K$ in a triangulated $PL$ manifold
$M$ is a closed embedded polygonal curve.
An (unoriented) {\em link} $L$ is a finite union of nonintersecting knots in $M$. 
Let $|L|$ denote the number of vertices of $L$. 

An {\em elementary move} (or $\Delta$-{\em move}) on a knot or link consists of one of the following operations, see
Burde and Zieschang \cite[p.~4]{BurZie85}, Murasugi \cite[p.~7]{Mur96}. Each of these operations is required to take
place in a single tetrahedron of $M$. 

(1) Split an edge into two edges by adding a new vertex to its interior. 

($1'$) [Reverse of (1)]
Combine two edges whose union is a line segment, erasing their common vertex.

(2) Given a line segment $AB$ of $L$ that contains a single interior vertex $D$ and a point $C$ not in $L$ such that
the triangle $ABC$ intersects $L$ only in edge $AB$, erase edge $AB$ and add two new edges $AC$ and $BC$ and new
vertex $C$. 

($2'$)
[Reverse of (2)]
If for two adjacent edges $AC$ and $BC$ of $K$ the triangle $ABC$ intersects $L$ only in $AC$
and $BC$ then erase
$AC$, $BC$ and vertex $C$ and add the edge $AB$, with a new interior vertex $D$.

\begin{figure}[hbtp]
\centering
\includegraphics[width=.5\textwidth]{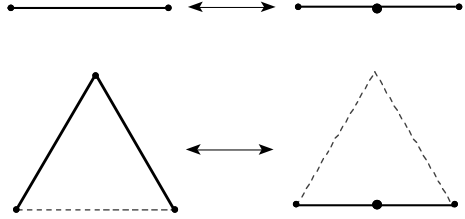}
\caption{
\label{elementary}
Two types of elementary moves.}
\end{figure} 

Corresponding to elementary moves of type (2), ($2'$) is the piecewise-linear operation of
linearly pulling the edges of a link across a triangle $ABC$. To define this operation, an extra
vertex $D$ on the segment
$AB$ must be specified; this is the role of elementary moves (1) and $(1')$. 

Two knots $K_1$ and $K_2$ are {\em equivalent} if one can be reached from the other by a finite
sequence of elementary moves.
This gives the same equivalence classes as requiring $K_1$ and $K_2$ to be {\em ambient isotopic} \cite{Rolfsen} in
$M$. The {\em unknot} is the class of knots equivalent to a triangle contained in a tetrahedron,
and a knot $K$ is {\em unknotted} if it is equivalent to the unknot. Equivalently, an unknot is
the boundary of a $PL$ embedded disk in $M$ \cite{Rolfsen}. We will work with knots embedded in
the interior of $M$. This is not a significant restriction, since it is possible to perturb any
knot which meets $\partial M$ into the interior of $M$ after some simplicial subdivision.

We conclude this section by outlining the proof of Theorem~\ref{thm102}. We first construct a compact
triangulated $PL$ submanifold $M_K$ by removing a solid torus neighborhood of $K$. By removing
the interior of a solid torus neighborhood $R_{K}$ of $K$ we construct a new triangulated
manifold $ M_K = M - R_K $. The boundary of $M_K$ consists of the original boundary $\partial M$
plus a 2-torus $\partial R_K$, called the {\em peripheral torus}.

If $K$ is unknotted then $M_K$ contains an {\em essential normal disk $S$}. Essential means that
$\partial S \subseteq \partial R_K$ is not the boundary of any disk in the torus $\partial R_K$.
Drawing on Haken's normal surface theory,
a result of Hass, Lagarias and Pippenger \cite{HasLagPip} implies that $S$
can be chosen to contain at most $ 2^{8t+6}$ triangles. 
We present this argument for completeness in \S2.
Its boundary $\partial S$ is a curve $K_2$. $K_2$ is
isotopic to $K$ within $\bar{R}_{K}$.
There is a sequence of at most $ 2^{8t+7}$
elementary moves which deform $K_2$ to a
triangle across triangles in $S$. This is carried out in \S3.
It remains to relate $K_2$ to $K$ by elementary moves. 

In \S3 we construct a curve $K_1$ isotopic to $K$ in $ R_K$ which is a {\em longitude},
a curve on the peripheral torus, $\partial R_K$, which bounds an embedded disk in $M_K$ but not inside the
peripheral torus.
Such curves exist if $K$ is unknotted in $M$.
In \S4 we construct a longitude
$K_2$ that lies in the 1-skeleton of the peripheral 2-torus $\partial R_K$,
consists of at most
$O(2^{c_3t})$ line segments and can be isotoped to $K$
using at most $O(2^{c_3t} )$ elementary moves.
In \S5 we prove results which are used to show that there is a
sequence of at most $O(2^{c_4 t})$ elementary moves which
deform $K_1$ to $K_2$ while remaining entirely
in the peripheral 2-torus $\partial R_K$.
In \S6 we combine these results to complete the proof.
All the constants and bounds above are explicitly computed.

The proof of Theorem~\ref{thm102} in effect deforms the original knot $K$ to
a trivial knot with
elementary moves made inside three surfaces:
an annulus is used to go from $K$ to the longitude
$K_1$, a 2-torus is used to go from $K_1$ to the
boundary of a normal disk $K_2$, and a disk is
used to go from $K_2$ to the trivial knot.

\section{Normal surfaces and elementary moves across an essential disk} \hsp

In this section we assume given a compact $PL$ 3-manifold with
(possibly empty) boundary ($M, \partial M$), 
which is triangulated using $t$ tetrahedra.
We will work with a knot $K$ embedded in the
1-skeleton of Interior$(M)$. By barycentrically subdividing $M$ twice, we obtain a manifold
$(M'', \partial M'')$ in which $K$ has a regular neighborhood $R_K$. The regular neighborhood, obtained by taking the
closed star of $K$ in $M''$, is a solid torus with
$\partial R_K$ disjoint from $\partial M$. We let $M_K = (M'' \backslash R_K ) \cup \partial
R_K$ and note that its boundary $\partial M_K = \partial M'' \cup
\partial R_K$. We call the manifold $(M_K , \partial M_K )$ with $\partial R_K$ marked a {\em truncated knot
complement} for $K$. 

A {\em meridian} on the peripheral torus $\partial R_K$ is an
oriented simple closed curve $\mu$ whose homology class
$[ \mu ]$ in $H_1 (\partial R, \ZZ )$ generates the kernel of the map
$$j_\ast : H_1 (\partial R_K , \ZZ ) \ra H_1 (R_K , \ZZ ),$$
where $j : \partial R_K \hookrightarrow R_K$
is the inclusion map. 
Equivalently, it is a simple closed curve which bounds a disk in $R_K$ but not in $\partial R_K$.

A {\em longitude} on the peripheral torus $\partial R_K$ is an
oriented simple closed curve $\lambda$ which intersects a
meridian transversely at a single point and 
whose homology class $[ \lambda ]$ in $H_1 (\partial R, \ZZ )$ generates the kernel of the
map
$$i_\ast : H_1 (\partial R, \ZZ ) \ra H_1 (M_K , \ZZ ),$$
where $i: \partial R_K \hookrightarrow M_K$ is the inclusion map, provided that
this kernel is nontrivial. Note that the kernel is always infinite cyclic or
trivial \cite{BurZie85}.

If the core $K$ of the solid torus $R_K$ is unknotted in $M$, then a longitude exists, and it bounds a disk in
$M_K$. For general $K$, a longitude does not exist when $K$ represents an element of infinite order in $H_1 (M , \ZZ
)$. If $M \cong S^3$, or more generally a homology sphere, then a longitude exists for any $K$.
This curve is sometimes called a {\em preferred longitude}.

When a longitude $\lambda$ of $\partial R_K$ exists, its
 intersection number on $\partial R_K$ with a meridian
$\mu$ is $\pm 1$, and the homology classes $[ \mu ]$
and $[ \lambda ]$ generate $H_1 (\partial R, \ZZ
) \cong \ZZ \ \oplus \ZZ$.

A {\em properly embedded disk} in a 3-manifold with boundary $(M, \partial M)$ is a disk $S$ PL-embedded in $M$
which satisfies $S \cap \partial M = \partial S$.
An {\em essential disk} in $M$ is a properly embedded disk such that
$\partial S$ does not bound a disk in $\partial M$, i.e.
the homotopy class of $\partial S$ in
$\pi_1 (\partial M)$ is non-trivial.

\begin{lemma} \label{knot}
Let $M_K$ be a triangulated $PL$ 3-manifold which is a truncated knot complement for $K$.
Then $K$ is unknotted if and only if $M_K$ contains an essential disk $S$ with $\partial S$ a longitude in $\partial R_K$. 
\end{lemma}
\bproof
If $K$ is unknotted then it bounds a disk which is embedded in $M$. This disk can
be intersected with $M_K$ to give an essential $PL$ disk $S$ in $M_K$. Conversely if $M_K$
contains an essential
$PL$ disk $S$ with $\partial S \subseteq \partial R_K $, then the boundary $\partial S$ of $S$
represents a longitude on $\partial R_K$, which is isotopic to $K$ in $R_K$ \cite[p.~29, Thorem
3.1]{BurZie85}, and therefore is equivalent to $K$. This disk can be used to contract
$\partial S$ until it becomes a single triangle in $M_K$.
The first homology group of a torus, $Z \oplus Z$, is canonically isomorphic to its fundamental
group, and the statement about the homology class of $[ \partial S ]$ is a consequence.$~~~\qed$ 

\begin{figure}[hbtp]
\centering
\includegraphics[width=.6\textwidth]{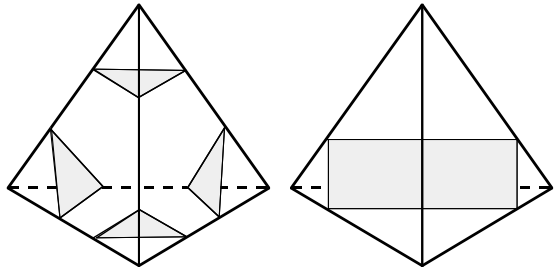}
\caption{
\label{normal}
Elementary disks in a normal surface.} 
\end{figure}

We will use normal surface theory in the form described in Jaco and Rubinstein 
\cite[Section~1]{JacRub89}.
A {\em normal surface} $S$ in
a triangulated compact 3-manifold $M$ is a $PL$-surface whose intersection with each tetrahedron in $M$ consists of
a disjoint set of {\em elementary disks}. These are either triangles and quadrilaterals. A quadrilateral consists of
two triangles glued together along an edge lying interior to the tetrahedron. All other edges of each triangle
or quadrilateral are contained in the 2-skeleton of the tetrahedron, and do not intersect any vertices. (See
Figure~\ref{normal}).

Within each tetrahedron of $M$ there are only 7 kinds of possible triangle or $PL$ quadrilateral, up to
$PL$-isotopies that map the tetrahedron to itself. These are labeled by the 7 possible ways of partitioning the
vertices of the tetrahedron into the two nonempty sets determined by cutting the tetrahedron along the triangle or quadrilateral. W. Haken observed that the isotopy-type of
a normal surface is completely determined by the number of
pieces of each of the 7 kinds that occur in each
tetrahedron\footnote{Haken's \cite{Hak61}
treatment of normal surfaces was based on a
formulation using cells. We follow an alternate development based on tetrahedra,
as in Jaco and Rubinstein \cite[Section~1]{JacRub89} and Jaco and Tollefson
\cite{JacTol95}}.
This can be described by a nonnegative integer vector
$\bv = \bv (S) \in \ZZ^{7t}$, which gives the
{\em normal coordinates} of $S$. He also showed that the set of
allowable values $\bv (S)$ for a normal surface $S$ lie in a certain homogeneous rational cone
${\sC}_M$ in $\RR^{7t}$ which we call the {\em Haken normal cone}.
If $\bv = ( v_1 , v_2 , \ldots
, v_{7t} ) \in \RR^{7t}$, then Haken's normal cone is specified by a set of linear equations and
inequalities of the form
\begin{equation} 
v_i \geq 0~, \quad \mbox{for}~~ 1 \leq i \leq 7t~, 
\label{eq401}
\end{equation} 
\begin{equation}
v_{i_1} + v_{i_2} = v_{i_3} + v_{i_4} \quad \mbox{(up to $6t$ equations)}~. 
 \label{eq402}
\end{equation}
The second set of equations expresses
{\em matching conditions} which say that the number of edges on a common triangular face of two adjacent tetrahedra coming from regions in each of the tetrahedra must match. For each triangular face there are three types of edges (specified by a pair of edges on the triangle), which yields 3~matching conditions per face. Triangular faces in the boundary $\partial M$ give no matching equations.
The cone ${\sC}_M$ is {\em rational} because all equations \ref{eq401}, \ref{eq402} have integer coefficients. We let
\begin{equation}
{\sC}_M ( \ZZ ) = {\sC}_M \cap \ZZ^{7t}~,
 \label{eq403}
\end{equation}
which is the set of integral vectors in the Haken normal cone ${\sC}_M$. Haken defined a {\em fundamental surface} to
be a normal surface $S$ such that 
\begin{displaymath} 
\bv (S) \neq \bv_1 + \bv_2 ~, \quad \mbox{with $\bv_1 , \bv_2 \in {\sC}_M ( \ZZ )
\backslash \{ \bf 0 \} $ } ~.
\end{displaymath} 
Such a vector $\bv (S)$ is called an element of the {\em minimal Hilbert basis}
${\sH}({\sC}_M )$ of ${\sC}_M$, using the terminology of integer programming, see Schrijver
\cite[Theorem~16.4]{Sch86}. 

\begin{theorem} \label{haken} (Haken)
If a triangulated, irreducible, orientable, compact 3-manifold with boundary
($M, \partial M$) contains an essential disk whose boundary is in a
component $\partial R$ of $\partial M$, then it contains such
an essential disk $S$ with boundary in $\partial R$
which is a fundamental normal surface.
\end{theorem}
\bproof
See Jaco-Rubinstein \cite[Theorem~2.3]{JacRub89} for a
treatment of the existence of an essential normal
disk in $M$ with boundary on $\partial R$.
For a fundamental disk, we refer to Jaco and Tollefson
\cite[Corollary~6.4]{JacTol95},
which establishes the stronger result that there exists such a disk that
is a fundamental surface and is also a {\it vertex}
surface, which means that $\bv (S)$ lies on an
extreme ray of Haken's normal cone ${\sC}_M$. $~~~\qed$

\paragraph{Remark.}
This result applies in arbitrary 3-manifolds.
The proof in \cite{JacTol95}
assumes compactness and irreducibility of $M$,
but these hypotheses can be removed. A classical
theorem of Alexander implies that knot complements in the 3-sphere are always irreducible.

Hass, Lagarias and Pippenger \cite{HasLagPip} give a simple bound for the complexity of any
fundamental surface, whose proof we include for completeness. 
\begin{lemma} \label{tetra}
Let $M$ be a triangulated compact 3-manifold, possibly with boundary, that contains $t$
tetrahedra.
\begin{itemize}
\item[$(1)$]
Any vertex minimal solution $\bv \in \ZZ^{7t}$ of the Haken normal cone $\hl{C}_M$ in $\RR^{7t}$ has 
\begin{equation}
\max_{1 \le i \le 7t} (v_i) \le 2^{7t-1} ~. 
 \label{eqN35}
\end{equation}
\item[$(2)$]
Any minimal Hilbert basis element $\bv \in \ZZ^{7t}$ of the Haken fundamental cone $\hl{C}_M$ has 
\begin{equation} 
\max_{1 \le i \le 7t} (v_i) < t \cdot 2^{7t+2} ~. 
\label{eqN36}
\end{equation}
\end{itemize}
\end{lemma}

\paragraph{Proof.}
(1)~Choose a maximal linearly independent subset of matching condition equations \ref{eq402}.
There will be $7t-d$ of them, where $d = \dim_{\RR} (\hl{C}_M )$.
Any vertex ray is determined
by adjoining to these equations $d-1$ other binding inequality constraints
\begin{displaymath} 
\{ v_{i_k} =0 ; ~~~1 \le k \le d-1 \} ~,
 \end{displaymath}
with the proviso that the resulting system have rank $7t-1$.
These conditions yield a $(7t-1) \times 7t$ integer
matrix $\hl{M}$ of rank $7t-1$, and the vertex ray elements $\bz = (z_1, \ldots, z_{7t} )$
satisfy
\begin{displaymath} 
\hl{M} \hl{z} = {\bf 0} ~.
\end{displaymath}
In order to get a feasible vertex ray in $\hl{C}_M$ all nonzero coordinates $z_i$ must have the same sign. 

At least one of the unit coordinate vectors $\be_1$, $\be_2 , \ldots \be_{7t}$ must be linearly independent of the
row space of $\hl{M}$. Adjoin it as a first row to $\hl{M}$ and we obtain a full rank $7t \times 7t$ integer matrix
$$\tilde{\hl{M}} = \left[ \begin{array}{c} {\be_k} \\ \hl{M}
\end{array} \right] , \quad \det ( \tilde{\hl{M}} ) \neq 0 ~. $$
Consider the adjoint matrix $adj ( \tilde{\hl{M}}) = \det ( \tilde{\hl{M}} ) \tilde{\hl{M}}^{-1}$, which has integer
entries 
\begin{equation}
w_{ij} = adj ( \tilde{\hl{M}})_{ij} = (-1)^{i+j} \det \tilde{\hl{M}} [j|i] ~, 
 \label{eqN39} 
\end{equation}
in which $\tilde{\hl{M}} [j|i]$ is the minor obtained by crossing out the $j^{\rm th}$ row and $i^{\rm th}$ column
of $\tilde{\hl{M}}$. 
Let 
\begin{displaymath} 
\bw = [w_{11}, w_{21}, \ldots, w_{7t,1} ]^t 
\end{displaymath}
be the first column of $adj ( \tilde{\hl{M}} )$. Since 
$ \tilde{\hl{M}} adj ( \tilde{\hl{M}} ) = \det ( \tilde{\hl{M}} ) \hl{I}$ this yields 
$$
 \hl{M} \bw = {\bf 0} ~,
$$
and $\bw \neq {\bf 0}$ because $\langle \be_k , \bw \rangle = \det ( \tilde{\hl{M}} ) \neq 0$.
We bound the entries of $\bw$ using {\em Hadamard's inequality}, which states for an $m \times
m$ real matrix $\hl{N}$ that 
\begin{equation}
\det ( \hl{N} )^2 \le
\prod_{i=1}^m \| \bn_i \|^2 ~,
\label{eqN311}
\end{equation}
in which $\| \bn_i \|^2$ is the Euclidean length of the $i^{\rm th}$ row $\bn_i$ of $\hl{N}$. We apply this to
equation~\ref{eqN39}, and observe that each row of the $(7t-1) \times (7t-1)$ matrix $ \tilde{\hl{M}} [j|i]$ has
squared Euclidean length at most 4, because this is true for all row vectors in the system of equations~\ref{eq402}
and for
$\be_k$. Applied to equation~\ref{eqN39} this gives
$$| w_{ij} |^2 \le 4^{7t-1} ~.$$
However equation~\ref{eqN39} shows that $\bw \in \ZZ^{7t}$, and a vertex minimal solution $\bv$ in the extreme ray
is obtained by dividing $\bw$ by the greatest common divisor of its elements, hence
$$\max_{1 \le i \le 7t} |v_i | \le \max_{1 \le i \le 7t} | w_{i1} | \le 2^{7t-1} ~. $$

(2)~A {\em simplicial cone} $\hl{C}$ in $\RR^{7t}$ is a $d$-dimensional pointed cone which has exactly $d$ extreme
rays.
Let $\bv^{(1)} , \ldots, \bv^{(d)} \in \ZZ^{7t}$ be the vertex minimal solutions for
the extreme rays. Each point in $\hl{C}$ can be expressed as a nonnegative linear combination of the $\bv^{(j)}$, as
$$ v = \sum_{j=1}^d c_j \bv^{(j)} , \quad \mbox{each} \quad c_j \ge 0 ~. $$
If $\bv$ is in the minimal
Hilbert basis ${\sH} (\hl{C} )$ of the cone $\hl{C}$, then $0 \le c_j \le 1$, for otherwise one has 
$$\bv = (\bv - \bv^{(j)} ) + \bv^{(j)} ~,$$ 
and both $\bv - \bv^{(j)}$ and $\bv^{(j)}$ are nonzero integer vectors in
$\hl{C}$, which is a contradiction. 
Thus, any minimal Hilbert basis element $\bv = (v_1, \ldots, v_{7t} )$ of a
simplicial cone satisfies
\begin{equation} 
|v_i | \le \sum_{j=1}^d | v^{(j)} | \quad \mbox{for} \quad 1 \le i \le 7t ~. \label{eqN312}
\end{equation}

The cone $\hl{C}_M$ may not be simplicial, but we can partition it into a set of simplicial cones $\{ \hl{C}_m^{(k)}
\}$ each of whose extreme rays are extreme rays of $\hl{C}_M$ itself. We have
$$
{\sH} (\hl{C}_M ) \subseteq \bigcup_k {\sH} (\hl{C}_M^{(k)} ) ~. $$
Thus all Hilbert basis elements of ${\sH} (\hl{C}_M )$ satisfy the bound \ref{eqN312} for Hilbert basis
elements of ${\sH} (\hl{C}_M^{(k)} )$. Using equation~\ref{eqN35} to bound $|v^{(j)}|$ we obtain 
$$|v_i | \le d2^{7t-1} < 7t 2^{7t-1} < t 2^{7t+2} ~,$$
as claimed.~~~$\qed$

\vspace*{+.1in}
A $PL$-embedded normal disk $S$ can be used as a template to transform its
boundary $\partial S$ to a single triangle by a
sequence of elementary moves, in such a way that all intermediate curves lie on the surface $S$.

\begin{lemma} \label{piece} 
Let $M$ be a piecewise-linear triangulated compact 3-manifold with boundary, and let $S$
be a normal disk in $M$ that consists of $w$ triangles. Then $\partial S$ can be isotoped to a triangle by a
sequence of at most $2w$ elementary moves on $ M$, each of which takes place in a triangle or edge
contained within $S$. \end{lemma}
\bproof
We will construct a sequence of $w$ simple closed curves $$
\partial S = C_0 , C_1 , C_2 , \ldots , C_{w-1} = \partial T $$
where $T$ is a triangle in $S$, and such that for $0 \leq i \leq w-1$:
\begin{enumerate}
\renewcommand{\labelenumi}{(\roman{enumi})} 
\item $C_{i+1}$ is obtained from $C_i$ by at most two elementary moves;
either a move of type $2'$ across a single triangle in $S$ followed by a move of type $1'$ which removes the
extra vertex created by the first move, or a move of type 1 followed by a move of type 2 across a single triangle in
$S$.
\item
There is a triangulated disk $S_i$ contained in $S$ which consists of $w-i$ triangles, such that $C_i = \partial
S_i$. \end{enumerate}

The construction proceeds by eliminating triangles in $S$ one at a time by a sequence of steps, each consisting
of two elementary moves, in such a way as to preserve the property that each $S_j$ is a
topological disk. We proceed by induction on $j$, with the base case
$j = 0$ holding by hypothesis. For the induction step, we note that the only way that removing from $S_j$ a triangle
that contains at least one edge of $\partial S_j$ can yield a surface $S_{j+1}$ that is not a topological disk is by
producing a vertex $\bv$ that is a cut-point, i.e. $\bv$ is visited twice in the curve $\partial S_{j+1}$. In that
case $\bv$ is a vertex of $\partial S_j$, which implies that the triangle that is removed must have exactly one edge
on $\partial S_j$, and the pair of elementary moves pulls this edge across to the two edges of this triangle that
are adjacent to $\bv$. We therefore proceed by first checking if $S_j$ has any triangle having two edges on
$\partial S_j$. If so, we remove this triangle to obtain $S_{j+1}$. If this cannot be done, then $S_j$
must contain an interior vertex, because a triangulated polygon that has no interior vertices and contains
more than one triangle contains at least 2 triangles having two edges on the boundary. In this case there exists a
triangle having one edge on the boundary which also has an interior vertex, and we can pull across this triangle to
obtain
$S_{j+1}$. This completes the induction step.
$~~~\qed$

Putting all these results together yields: 
\begin{lemma} \label{disk.count}
Let $M_K$ be a triangulated $PL$ 3-manifold consisting of $t$
tetrahedra which is a truncated knot complement
for $K$. If $K$ is unknotted, then there is an essential normal disk $S$ in $M_K$ that contains
at most $ 2^{8t+6}$ triangles, and $S$ has a
polygonal boundary $ K_1 = \partial S$ which lies on the peripheral
2-torus $\partial R_K$ and consists of at most $ 2^{8t+7} $
line segments.  The boundary $K_1$ of $S$ can be
transformed to a single triangle by at most $ 2^{8t+7}$ elementary moves,
with all intermediate curves being the boundary of some triangulated subdisk of $S$.
\end{lemma}
\bproof By Lemma~\ref{knot} and Theorem~\ref{haken} there exists an essential normal disk which is a fundamental
surface for $M_K$. 
The bound of Lemma~\ref{tetra} implies that this surface contains at most $7t^2 \cdot 2^{7t+2} 
< t^2 \cdot 2^{7t+5} < 2^{8t+5}  $ triangles and quadrilaterals. Since each quadrilateral consists of two triangles in the
triangulation, we obtain the upper bound $2^{8t+6} $
for the total number of triangles. The boundary $\partial S$
contains at most two edges of each triangle,
for a total of at most $2^{8t+7}$ edges.
Lemma~\ref{piece} gives
the stated bound on elementary moves. $~~~\qed$

\paragraph{Remark.}
We could improve the above estimate by using a vertex minimal solution which is a disk.

\section{Longitudes and cores of solid tori} \hsp
In this section we bound the number of elementary moves needed to isotop a knot $K$ in the 1-skeleton of a
triangulated compact 3-manifold $(M, \partial M)$ across an annulus
to a curve $K_2$ which is a longitude on the peripheral torus
$\partial R_K$ of $K$. The main result of this section is the following: 

\begin{theorem}
\label{longitude.count} 
Let $(M, \partial M)$ be a triangulated compact $PL$ 3-manifold with $t$ tetrahedra, let
$K$ be an unknotted curve embedded in Interior$(M)$,
and let $\partial R_K$ be the peripheral torus of $K$ in
the second barycentric subdivision of $M$.
Then there exists a $PL$ closed curve $K_2$ in $\partial R_K$ which is a
longitude for $\partial R_K$ in $M_K$ such that: 
\begin{description}
\item{\hspace*{.25in}(i).}
There is an embedded $PL$ annulus $S$ in the solid torus $R_K$ whose two boundary components are
$K_2$ and $K$, and which consists of at most $2^{1858t}$ triangles. 
\item{\hspace*{.25in}(ii).} There is
an isotopy from $K_2$ to $K$ across $S$ that consists of at most
$2^{1858t+1}$ elementary moves.
\end{description}
\end{theorem}
\paragraph{Remarks.}
The closed curve $K_2$ is embedded in $\partial R_K$ but generally does not lie in its 1-skeleton.
In the special case that $(M, \partial M)$ is $PL$-homeomorphic to a triangulated convex polyhedron in $\RR^3$, as
occurs with a standard knot, the exponential bounds in (i), (ii) can be improved to bounds polynomial in $t$, as we
indicate after the proof. 

The proof of Theorem~\ref{longitude.count} proceeds in three steps.
We first construct a closed curve $\alpha$ embedded in the 1-skeleton of $R_K$ which is {\em parallel} to the core
$K$, i.e. its homology class $[ \alpha ]$ in $H_1 (\partial R_K , \ZZ )$ has intersection number $\pm 1$ with a
meridian $[\mu ]$. We next observe that the
homology class in $H_1 (\partial R_K , \ZZ )$ of a longitude $K_2$ is necessarily
\begin{equation} 
[K_2] = \pm ([ \alpha ] + k [ \mu ])
\label{eq401a}
\end{equation}
for some integer $k$.
The second step is to bound $|k|$ by an exponential function of $t$. The third step is to obtain $K_2$ from
$\alpha$ by adding $k$ twists to $\alpha$ near a fixed meridian $\mu$. 
Given $K_2$, we construct a triangulated
annular surface $S$ embedded in $R_K$ which lies in the 2-skeleton of $R_K$ except in a ``collar'' of one meridian disk,
into which the $k$ twists are inserted. All vertices of the triangulated annulus $S$ lie on its boundary $\partial
S = K_2
\cup K$. Finally we obtain an isotopy giving the elementary move bound by pulling $K$ across the triangles in $S$ to
$K_2$. 

We begin with some preliminary lemmas involving $PL$ topology, which pin
down some of the combinatorial structure of the regular neighborhood of $K$.
Background for $PL$ topology can be found in \cite{Hudson} or \cite{{RS}}. 

Let $v_0 , \ldots , v_r$ be the vertices of the knot $K$ in $M$. Since $R_K$ lies in
the second barycentric subdivision of $M$, each edge $[v_i , v_{i+1} ]$ of $K$ is subdivided into four edges by the
addition of three new vertices . We will continue to denote the subdivided
curve by $K$. Denote the original vertices
$v_i$ by $w_{4i}$ and let $w_{4i+1}, w_{4i+2}, w_{4i+3}$, $0 \le i \le r$ denote the new vertices added in the two barycentric
subdivisions. We also use the convention that $w_{s} = w_0$, where $s = 4r+4$ is the total number of vertices, so that the vertices
$w_i$ of $K$ are cyclically ordered.

 We will break up the solid torus $R_K$ into a union of what we call {\em wheels} and {\em wedges}. For each edge
$[w_{i-1} , w_i ]$ let the {\em wheel} around that edge be the closed star $St( [w_{i-1} , w_i ], R_K )$ of the edge 
$[w_{i-1} , w_i ]$ in the solid torus $R_K$. This is the union of all the tetrahedra in $R_K$ having $[w_{i-1} , w_i ]$ as
an edge, and is homeomorphic to a closed ball. For each vertex $w_i$ on $K$ let the {\em wedge around $w_i$}
be the union of all the tetrahedra in $R_K$ meeting $K$ only at the point $w_i$.
Let
\begin{equation} 
\mu_i= lk ( [w_{i-1} , w_i ] , R_K )~, \quad 1 \leq k \leq s~,
\end{equation}
be the link of the edge $[w_{i-1} , w_i ]$ in the solid torus $R_K$.

\begin{figure}[hbtp]
\centering
\includegraphics[width=.4\textwidth]{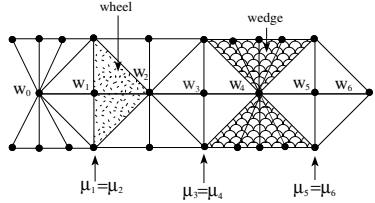}
\caption{A 2-dimensional section of $R_K$.} 
\label{section} 
\end{figure}

\begin{lemma} \label{lem.meridian}
(1) Each $\mu_i$ lies in $\partial R_K$ and is a meridian. \\
(2) $\mu_i$ and $\mu_{j}$ coincide if $i$ is odd and $j=i+i$ or if $i$ is even and $j=i-1$. Otherwise $\mu_i$ and $\mu_j$ are
disjoint. \\ 
(3) An edge on $\partial R_K$ either lies on some $\mu_i$ or else lies on a unique triangle in $R_K$ whose
opposite vertex $w$ lies on $K$.
\end{lemma}

\begin{figure}[hbtp]
\centering
\includegraphics[width=.6\textwidth]{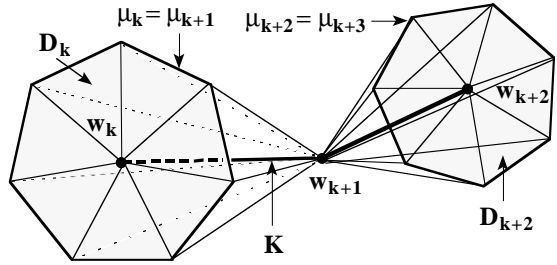}
\caption{
\label{meridian} Adjacent meridians. The index $k$ is odd.} 
\end{figure}

\bproof
(1) The closed star $St( [w_{i-1} , w_i ], R_K )$ is supported on the set of tetrahedra in $R_K$ having 
$[w_{i-1} , w_i ]$ as an edge. Each such tetrahedron has one edge which is disjoint from 
$[w_{i-1} , w_i ]$, and the link $\mu_i$ is the union of these edges. The surface $D_i$ which is obtained by
coning $\mu_i$ to the vertex $w_i$ forms a topological disk in
$R_K$ whose boundary is $\mu_i$, and which intersect
$K$ transversely in a single point. Thus $\mu_i$ is a meridian; see Figure~\ref{meridian}.

(2) This follows since $R_K$ is a barycentric subdivision. The vertices $w_i$ with $i$ odd were added during the second barycentric
subdivision. Prior to the second subdivision, these vertices lay along interior points of edges of $K$. The two meridians associated to the
edges $ [w_{i-1} , w_i ]$ and $ [w_{i} , w_{i+1} ]$ coincide;
see Figure~\ref{section} for a cross sectional sketch of $R_K$.
For $i$ even, the distance between $w_{i-1}$ and $w_{i+1}$ in $link(w_i)$ is at least four, since the link
has been barycentrically subdivided. If
$\mu_i \cap \mu_{i+1} \neq \emptyset$ then this distance equals two.

(3) Each edge in $\partial R_K$ is an edge of a triangle with a vertex $w_i$ on $K$, since $ R_K$ is obtained by taking the
closed star of $K$ in $M''$. If the edge is in a meridian $\mu_i$, then it lies in triangles with vertices at
both $w_{i-1}$ and $w_i$ or $w_i$ and $w_{i+1}$. If it doesn't lie on a meridian 
then it is contained in the wedge around $w_i$. Any triangle containing the edge and a point $w'$ on $K$ cannot
cross the two wheels meeting $w_i$, and the point $w'$ must coincide with $w$.
$~~~\qed$

\begin{lemma} \label{parallel curve}
There exists an embedded closed curve $\alpha$ in the 1-skeleton of $\partial R_K$ which is parallel to $K$ in $R_K$.
There is an annulus $S_0$ embedded in the 2-skeleton of the solid torus $R_K$ whose boundary components
are $\alpha$ and $K$, and all vertices of $S_0$ lie in its boundary $\alpha \cup K$. 
\end{lemma}
\bproof
 We construct $\alpha$ as the union of a sequence of arcs. For each odd $i, \ 1 \le i \le s$, we pick a path $\alpha_i$ in
the 1-skeleton of $\partial R_K$ that connects a vertex $x_i$ of $\mu_i = \mu_{i+1}$ to a vertex $x_{i+1}$
on $\mu_{i+2} = \mu_{i+3}$, does not cross $\mu_{i-1}$, and is shortest among all such paths. (We don't allow $\alpha_i$
to cross $\mu_{i-1}$ to prevent it
from running between the curves in the wrong direction.) 
We set $x_s = x_0$ as before. Note that
$\alpha_i$ is embedded and that it meets the collection of meridians only at its endpoints, since otherwise it can be shortened. 

We then add for each odd $i$ a connecting path $\beta_i$ on $\mu_i$ from the endpoint $x_{i-1}$ of
$\alpha_{i-2}$ to the initial point $x_i$ of $\alpha_i$, and we set
\begin{equation} 
\alpha = \beta_1 \alpha_1 \beta_3 \alpha_3 \beta_5 \ldots \beta_{s-1}\alpha_{s-1} ~. \end{equation}
The path $\alpha$ is a closed curve, and is embedded in $\partial R_K$.

We construct the surface $S_0$ as follows. Each path $\beta_i$, $1 \le i \le s$, $i$ odd, lies
on a meridian $\mu_i = \mu_{i+1}$ which
is the boundary of a meridian disk $D_{i+1}$ centered at the vertex $w_{i+1}$ on $K$, and we cone $\beta_i$ to $w_{i+1}$ along
$D_{i+1}$, obtaining a subdisk of $D_{i+1}$. For each edge $\alpha_i$ that lies on the boundary of the wedge around the vertex
$w_{i+1}$, there is a unique triangle containing that edge and
$w_{i+1}$. For fixed $i$, the union of these triangles forms a disk in $R_K$ meeting $K$ at the point $w_{i+1}$, and $\partial R_K$
along $\beta_i \cup \alpha_i$, and whose boundary consists of $\beta_i \cup \alpha_i$ together with two interior edges of $R_K$, one
running from the initial point $x_{i-1}$ of $\beta_i$ to
$w_{i+1}$ and one running from the endpoint $x_{i+1}$ of $\alpha_{i}$ to $w_{i+1}$ .

\begin{figure}[hbtp]
\centering
\includegraphics[width=.6\textwidth]{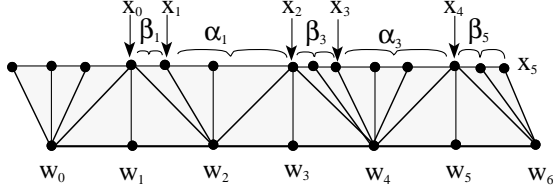}
\caption{
\label{annulus} The annulus $S_0$.} 
\end{figure}

Finally for each meridian $\mu_i$ with $i$ odd,
there is a pair of triangles
$$[w_{i-1}, w_i, x_{i-1}]$$
and
$$[ w_i, w_{i+1},x_{i-1}],$$
whose union is a disk contained in the two wheels around $[w_{i-1}, w_i]$ and $[ w_i, w_{i+1}]$, and connecting the free edge running
from the endpoint $x_{i-1}$ of $\alpha_{i-2}$ to $w_{i-1}$ to the free
 edge running from the
initial point (also $x_{i-1}$) of $\beta_i$ to $w_{i+1}$. 
The union of all these disks forms an annulus which is embedded in the 2-skeleton of
$R_K$, whose boundary components are $\alpha$ and $K$, and all vertices of which lie in $\alpha \cup K$. 
See Figure~\ref{annulus}. $~~~\qed$ \\

{\bf Proof of Theorem~\ref{longitude.count}.}
We wish to convert the closed curve $\alpha$ given in Lemma~\ref{parallel curve} to a longitude $K_2$ by adding
``twists''. The homology class $[ K_2 ]$ of a longitude $K_2$ satisfies
\begin{displaymath} 
[ K_2 ] = \pm ([ \alpha ] +k [ \mu ]) 
\end{displaymath}
where $\mu$ is a fixed meridian.
According to Lemma~\ref{lem.meridian} we may choose $\mu$ to be $\mu_1$ taken with a fixed orientation.

The number $k$ represents the number of ``twists'' that we must insert in $\alpha$ to convert it to the knot $K_2$.
We proceed to bound $k$ in terms of the number of tetrahedra $t$ in the original manifold $M$. We may assume $k
\neq 0$ for otherwise we are done. The value of $k$ is uniquely determined by the requirement that the homology
class $[ \alpha ]+k [ \mu ]$ vanish in $H_1 (M_K , \ZZ )$. (Here we use the fact that $[ \mu ]$ is nonvanishing
in $H_1 (M_K , \ZZ )$.) This is equivalent to the curve $\alpha + k \mu$ being a boundary in the cycle group.
The manifold $M_K$ contains at most $576t$ tetrahedra, since it is contained in the second barycentric subdivision,
hence the number $e$ of edges in $M_K$ is at most $3456t$.
The group of boundaries in the $PL$ chain complex of $M_K$
is generated by the paths around all 2-simplices in $M_K$, and the number $f$ of these is at most
$2304t$. We choose
fixed orientations for the edges and faces, and we can then describe a basis of generators for the group of
boundaries using an $e \times f$ matrix $\hl{M} = [ \hl{M}_{ab} ]$ with entries
$$
\hl{M}_{ab} = \left\{
\begin{array}{ll}
\pm 1 & \mbox{if $a$ is an edge in the oriented 2-simplex $b$}~, \\ 0 & \mbox{otherwise}~.
\end{array}
\right.
$$
in which each column of $\hl{M}$ has exactly three nonzero entries. Let $\bc_\alpha$ be the $e \times 1$-vector describing $\alpha$
as an oriented cycle in $M_K$, and $\bc_\mu$ be the $e \times 1$ vector describing
$\mu$ as an oriented cycle in $M_K$. All entries in $\bc_\alpha$ and $\bc_\mu$ are 0,1 or $-1$. The condition for
$\bc_\alpha + k \bc_\mu$ to be a boundary is that the linear system
\begin{equation}
\hl{M} \bx = \bc_\alpha +k \bc_\mu
\label{eq407a}
\end{equation}
have a solution in integer unknowns $\bx = (x_1 , \ldots , x_f )^T$. The hypotheses imply that this linear system is
solvable over $\ZZ$ for a unique value of $k$. By taking a generating set of boundary faces, we have that $\hl{M}$
is of full column rank. The meridian class
$\bc_\mu$ is not in the column space of $\hl{M}$ since otherwise $\mu$ would bound a surface in $M_K$ and the
unknotted curve $K$ would represent an element of infinite order in $H_1 (M , \ZZ )$. It
follows that there is at most one {\em real} value of
$k$ for which the system (\ref{eq407a}) is solvable, for if there were two
then $\bc_\mu$ would be in the column space of $\hl{M}$.
Thus $\bc_\mu$ is linearly independent of the column space
of $\hl{M}$ and the linear system
$$
[ \hl{M} , \bc_\mu ] \left[ \df{\bx}{-k} \right] = \bc_\alpha ~, $$
has full column rank.
Choosing a suitable $f+1$ rows of
$[ \hl{M} \bc_\mu ]$ yields an invertible matrix $\tilde{\hl{M}}$ and 
\begin{displaymath} 
[ \df{x}{-k} ] = \tilde{\hl{M}}^{-1} \tilde{\bc}_\alpha ~, 
\end{displaymath}
in which $\tilde{\bc}_\alpha$ is the subset of $\bc_\alpha$ determined by the same $f+1$ rows.
By Cramer's rule we have
\begin{displaymath} 
-k = \df{1}{\det (\tilde{\hl{M}} )}
\sum_{i=1}^{f+1} \tilde{\hl{M}} [i,f+1] (\tilde{\bc}_\alpha )_i ~, 
\end{displaymath}
where $\tilde{\hl{M}} [i, f+1]$ denotes a minor of $\tilde{\hl{M}}$ on the last column.
Since $| (\tilde{\bc}_\alpha )_i | \leq 1$ we have 
\begin{equation}\label{eq410a}
|k| \leq \sum_{i=1}^{f+1} | \tilde{\hl{M}} [i, f+1 ] |~. \end{equation}
We bound the terms in this expression using Hadamard's determinant bound (\ref{eqN311}).
Each column of $\tilde{\hl{M}}$ has at
most three nonzero entries, each equal to $\pm 1$, except the last column, which does not appear in the $f \times f$ minor $\tilde{\hl{M}} [i, f+1 ]$.
This yields
the bound
\begin{displaymath} 
| \tilde{\hl{M}} [i, f+1] |^2 \leq 3^f ~, ~~1 \le i \le f+1 ~,
\end{displaymath}
which with (\ref{eq410a}) yields
\begin{equation}
|k| \leq (f+1) 3^{\frac{f}{2}} \leq (2304t+1) 3^{1152t}~. 
\label{eq412a}
\end{equation}

We now describe how the ``twists'' are added to $\alpha$ and the surface $S_0$ of Lemma~\ref{parallel curve}.
to yield an annulus $S$ with boundary a longitude.
The disk $D_2$ cuts the core $K$ at the
point $w_2$ and stretches to the boundary. We add the twists to $S_0$ by modifying it inside the wedge around
$w_2$. 
Let $A$ denote the surface consisting of the set of closed triangles in the wedge around
$w_2$ which intersect the meridian $\mu_2$ and which lie in
$\partial R_K$. Then $A$ is an annulus, and we suppose it contains $\ell$ triangles.
The first edge $e$ of $\alpha_1$ has one endpoint on $\mu_1 = \mu_2$ and the other endpoint disjoint from
$\mu_2$ and on the other boundary component of $A$. 

\begin{figure}[hbtp]
\centering
\includegraphics[width=.5\textwidth]{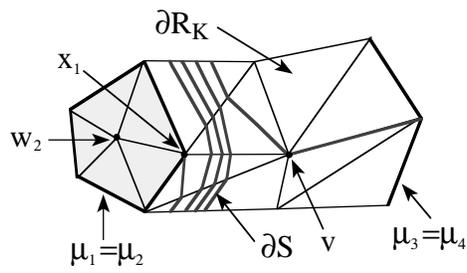}
\caption{
\label{annulus2} Adding twists replaces $S_0$ by $S$.} 
\end{figure}

We replace the edge $[x_1, v]$ of $\alpha$ with a path $\alpha'$ of $k \ell $ edges which form a spiral going $k$ times around
$A$, starting at $x_1$ and ending at $v$; see Figure~\ref{annulus2}.
To convert the surface $S_0$ of Lemma~\ref{parallel curve} to the desired surface $S$, we first subdivide the interval
$[x_1 ,v ]$ into $k \ell $ segments by adding vertices 
$y_1 , y_2 , \ldots , y_{k \ell - 1}$.
We remove the triangle $[x_1 , v , w_2]$ from $S_0$ and add $k \ell $ new triangles,
obtained by coning the edges of $\alpha'$ to $w_2$. This construction produces
the surface $S$, which is also an annulus and which also has
its vertices either on $K$ or $\partial R_K$. 

\begin{figure}[hbtp]
\centering
\includegraphics[width=.5\textwidth]{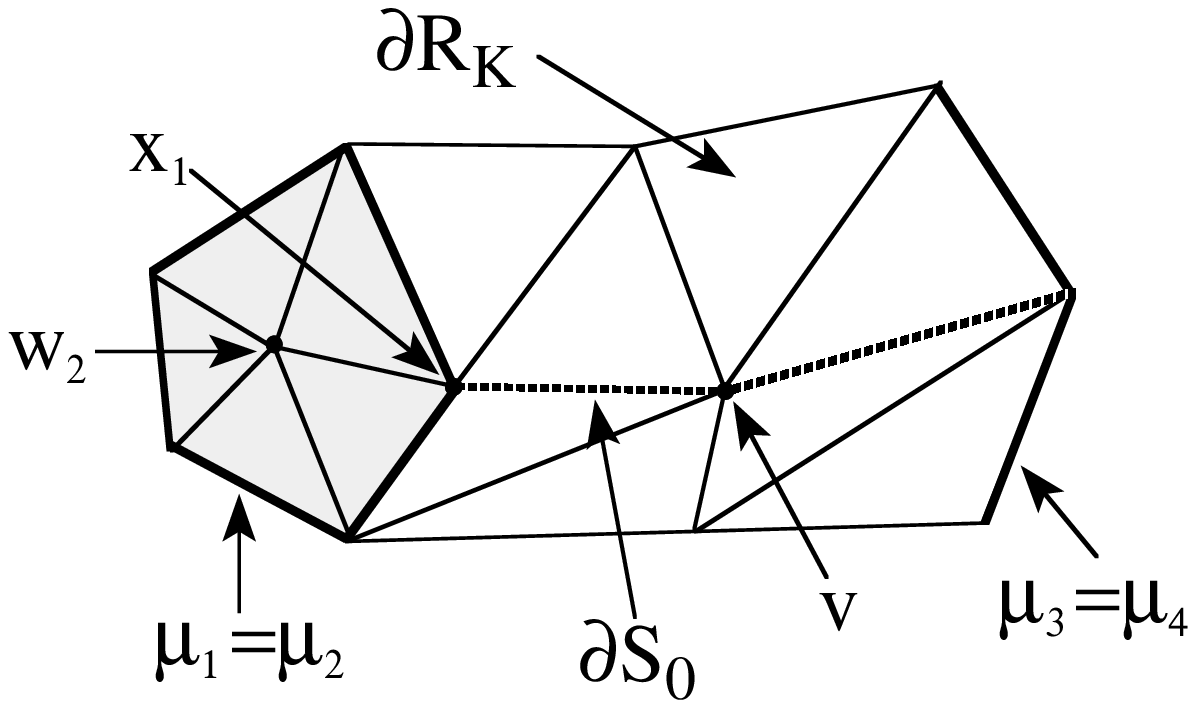}
\caption{
\label{annulus1} $S_0$.} 
\end{figure}
 
It is now straightforward to complete the proof. To prove (i) we observe that the number of
triangles in $S$ is at most the number of triangles in $R_K$, which is at most $2304t$, plus the
number $k \ell$ of extra triangles added in the ``twisting'' step in constructing $S$ from
$S_0$. The trivial bound for $\ell$ is the number of triangles in
$M''$ which is at most $2304t$. Combining this with (\ref{eq412a}) yields at most

\begin{displaymath} 
2304t + k \ell \le 2304t + (2304t + 1) 3^{1152t} (2304t) 
\leq t^2 \cdot 3^{1152t + 20} < 2^{1858t}
\end{displaymath}
triangles in $S$.

To prove (ii) we pull the curve $K$ across $S$ to $K_2$ one triangle at a time.
Using the fact that all vertices of $S$ are in the boundary $\partial S$, a suitable legal order of triangles to pull across is easily found.
Counting two elementary moves per triangle, this yields at most $2^ {1858t+1}$
elementary moves in total.
$~~~\qed$ \\ 
\paragraph{Remark.}
In the special case that the manifold $M$ is embedded in $\RR^3$, as when we start with a
standard knot complement, we can improve the bounds (i), (ii) of Theorem~\ref{longitude.count} to be polynomial in $n$. We
sketch the idea. The knot $K$ has at most
$2304t$ edges in $M$. Using the embedding of the manifold $M$ in $\RR^3$ we can find a regular projection of the
knot to a knot diagram
$\sD$, and this has
$O(n^2 )$-crossings, as we show in Lemma~\ref{regular.proj} of \S7.
We can specify a curve $\alpha$ in the 1-skeleton of $\partial R_K$ which is parallel to the core $K$ whose projection runs just to the right of $K$ in the projection plane.
Choose an orientation of $K$.
A longitude of a knot $K$ in $\RR^3$ is characterized by having a linking number with $K$ equal to zero. The
linking number is easily calculated from the projection. To each crossing of $K$ is associated a sign, and the sum of
these signs, called the {\em writhe} of the knot, gives the linking number of a curve with $K$. To ensure that the linking
number is zero, we alter the construction of $\alpha$.
Near every overcrossing of $K$, we add a full twist of $\alpha$
around $K$.
This twist is positive if the sign of the crossing is negative, and vice versa.
The result is a knot $K_2$
in $\RR^3$ that lies on the peripheral torus whose linking number with $K$ is zero.
We can construct the associated annulus $S$ in $R_K$ with boundary
components $K$ and $K_2$ similar to the proof given for Theorem~\ref{longitude.count}, and it has $O(n^2 )$ triangles.
In this argument the key step is that the writhe is bounded in terms of the number of crossings of $\sD$, which is $O(n^2
)$, and this gives a bound $k = O(n^2 )$ in place of (\ref{eq412a}).

\section{Isotopies of curves on a surface} \hsp 
In this section we examine isotopies between two $PL$ curves $\alpha$ and
$\beta$ on an oriented triangulated surface $F$, and bound the number of elementary moves needed to move $\alpha$ to
$\beta$. In the application,
the surface will be the peripheral torus $\partial R_K$ of a knot $K$,
$\alpha$ will be $K_1$,
the boundary of a normal disk, isotopic on $\partial R_K$ to $\beta$, and
$\beta$ will be a fixed longitudinal curve $K_2$,
parallel to $K$ in its regular neighborhood $R_K$.
We derive the results of this section for a general surface,
since no extra work is needed. 

Our main result is the following.

\begin{theorem} \label{surface.bound.thm} 
Let $F$ be a triangulated orientable surface with $u$ triangles, and let $\alpha$ and $\beta$ be
embedded isotopic $PL$ curves on $F$, each having at most $s$ line segments.
Then there is an isotopy carrying $\beta$ to $\alpha$ via no more than
$ 2^{11} (s + u^2)^4 u^3$
elementary moves.
\end{theorem}

\begin{figure}[hbtp]
\centering
\includegraphics[width=.3\textwidth]{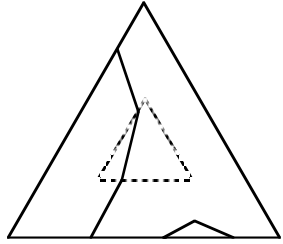}
\caption{
\label{basic.arc} Type 1 and type 2 basic arcs.} 
\end{figure}

\begin{corollary} \label{surface.bound.cor} 
Let $F$ be a triangulated orientable surface with at most
$2^{10}t$ triangles, and let $\alpha$ and $\beta$ be
embedded isotopic $PL$ curves on $F$, each having at most $2^{1070209t}$
line segments.
Then there is an isotopy carrying $\beta$ to $\alpha$ via no
more than $2^{10^7t -1}$ elementary moves.

\end{corollary}
We will work with specially positioned curves and isotopies, called {\em basic}, which we now define.
Let $E$ be a unit edge length equilateral triangle. Let $E'$ be an equilateral triangle of
half the size; see Figure~\ref{basic.arc}.
A {\em basic arc}
in $E$ is one of two types:

(a)~{\em Type 1.}
An arc running between the interiors of two distinct edges of $E$. This type of arc contains three straight line
segments. The first segment runs from the initial point towards the center of $E$ until it meets
$E'$. The third segment runs from the final point of the arc towards the center of $E$ until it
meets $E'$. The middle line segment connects the two in $E'$. \\

(b)~{\em Type 2.}
An arc with both of its endpoints on the interior of a single edge of $E$.
This type of arc
contains two straight line segments, each of which makes an angle of $1^o$ with the edge. \\

Arcs of type 2 are disjoint from any arcs of type
1 with non-linking endpoints, since an arc of type one makes an angle larger than $30^0$ with
an edge of $E$, and arcs of type 2 are disjoint from the smaller equilateral triangle.

Now let $E'$ be any triangle and fix a linear map from $E$ to $E'$.
A {\em basic arc in $E'$} is defined to be the image of a basic arc in $E$.

A {\em basic curve} in a triangulated surface is a curve which meets each triangle in a union of
basic arcs.

\begin{lemma} \label{basic.arcs}
Any two basic arcs in a triangle with distinct endpoints which lie on the
interiors of the edges of $E$ intersect in either zero or one point.
Any collection of disjoint embedded curves in a triangle can be isotoped, without moving their
endpoints, to a collection of embedded basic arcs. 
\end{lemma}
\bproof Two basic arcs with distinct endpoints intersect if and only if their endpoints link
on the boundary of the triangle. Any family of disjoint embedded arcs in a triangle is isotopic,
rel boundary, to any other such family. $~~~\qed$ 

\begin{figure}[hbtp]
\centering
\includegraphics[width=.5\textwidth]{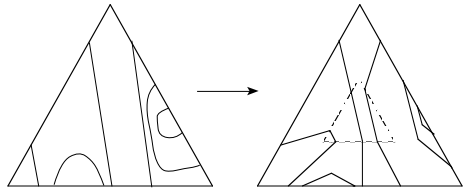}
\caption{
\label{basic.curves} Any set of arcs can be isotoped to a set of basic arcs.} 
\end{figure}
 
A {\em basic move} of a basic curve in a triangulated surface is defined to be one of the
following curve isotopies, each of which produces a new basic curve:

\begin{figure}[hbtp]
\centering
\includegraphics[width=.5\textwidth]{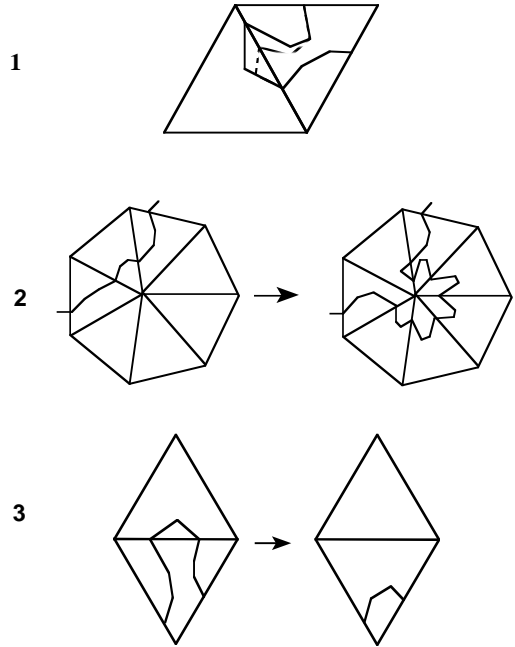}
\label{basic.moves}
\caption{ Basic move of types 1, 2 and 3.}
\end{figure}

(a)~{\em Type 1.}
An isotopy that slides an arc of the curve that crosses an edge along that edge.

(b)~{\em Type 2.}
An isotopy that
isotops an arc of the curve that crosses an edge of the triangulation across an adjacent vertex, and then isotops each of
the arcs of intersection with a triangle created or altered in this process to basic arcs.

(c)~{\em Type 3.}
An isotopy that slides an arc of the curve with both endpoints on one edge of a triangle $E$ across the 2-gon it forms
with that edge. This results in the joining of two arcs in an adjacent triangle $E'$. If the two arcs give a
closed curve in $E'$, the closed curve is removed. Otherwise, this new arc is
isotoped in $E'$ to a basic arc.

We count the number of elementary moves required by each basic move.
We first look in a single triangle.

\begin{lemma} \label{polygonal.delta}
Let $\alpha =\alpha_1 \cup \alpha_2 \cup... \cup \alpha_n $ and $\beta =\beta_1 \cup \beta_2 \cup... \cup \beta_n$
each be a union of $n$ connected polygonal arcs in a triangle
with the endpoints of $\alpha_i$ and
$\beta_i$ equal, and having a total of $a$ and $b$ edges respectively. Then
$\alpha$ can be isotoped to $\beta$ in the triangle using at most $a+b$ type 2 elementary moves, and $2(a+b)$
elementary moves in total.
\end{lemma}
\bproof 
We will count the number of required elementary moves in two cases. In each case we will isotop both 
$\alpha$ and $\beta$ to an intermediate curve $\gamma$, which stays close to the boundary of the triangle. To construct
$\gamma$, we first pick $\theta$ to be a positive angle which is small enough so that for any vertex $v$ in $\alpha \cup
\beta$ which lies on an edge of the triangle and any vertex $w$ in $\alpha \cup \beta$ which does not lie on that edge,
the angle between that edge and the line segment $vw$ is greater than $\theta$.

Given an arc $\alpha_i$ of $\alpha$, we construct an arc $\gamma_i$ with the same endpoints, but with interior disjoint
from $\alpha$ and containing just two segments. If $\alpha_i$ runs between distinct edges of the triangle
then it cuts off a subdisk containing a single vertex of the triangle. Construct $\gamma_i$ to consist of two
segments, each of which makes an angle of $\theta$ with the edges going towards that vertex. By definition of $\theta$,
$\gamma_i$ will be disjoint from the interiors of $\alpha_i$ and $\beta_i$. If $\alpha_i$ has both endpoints on a single
edge of the triangle, then it cuts off a subdisk containing no vertices and a segment of one edge. Construct $\gamma_i$ to
consist of two segments contained in this disk, each of which makes an angle of $\theta$ with its edge. Again,
$\gamma_i$ will be disjoint from the interiors of $\alpha_i$ and $\beta_i$.

\begin{enumerate}
\item $\alpha$ and $\beta$ are connected.

If $\alpha$ does not meet both sides of $\beta$, then the region between them can be triangulated with $a+b-2$ triangles,
which can be used to isotop $\alpha$ to $\beta$ using $a+b-2$ type 2 elementary moves. If they intersect, then
applying the previous argument shows that each can be isotoped to the disjoint curve $\gamma$ by using $a$ and $b$
type 2 elementary moves respectively. It follows that we can always isotop $\alpha$ to $\beta$ using at most $a+b$
type 2 elementary moves, and $2(a+b)$ elementary moves of all types.

\item $\alpha$ and $\beta$ are not connected.
 
We will construct an isotopy from $\alpha$ to $\beta$ which keeps the curves embedded at each stage.
Let $a_i$ be the number of segments of $\alpha_i$, so that $\sum_{i=1}^n a_i = a$. 
Each arc $\alpha_i$ of $\alpha$ cuts off a disk meeting exactly one vertex of the triangle in which it lies, if it runs
between distinct edges of the triangle, or no vertices if it runs from an edge back to that edge.
An arc of $\alpha$ is
said to be {\em outermost} if it cuts off a disk which contains no other arcs of $\alpha$.
Let $\alpha_i$ be outermost. Then $\gamma_i$ is disjoint from all the other arcs in $\alpha \cup \beta$
and meets $\alpha_i \cup \beta_i$ only in their endpoints. The number of segments of $\gamma_i$ is two,
and $\gamma = \bigcup \gamma_i$ is embedded. 
We can isotop the outermost arc $\alpha_i$ of $\alpha$ to the 
arc $\gamma_i$ with the same endpoints without ever introducing any
new intersections in $\alpha$,
since the other arcs in
$\alpha$ lie on the other side of $\alpha_i$ from $\gamma_i$.
This isotopy requires at most $a_i$ type 2 elementary
moves, as in the previous case. 
We then repeat for an outermost arc $\alpha_j$ among the collection of arcs of
$\alpha - \{\alpha_i\}$. This isotopy takes place between $\alpha_j$ and $\gamma_j$,
and in particular never meets
$\gamma_i$. Repeating for all arcs in $\alpha$, we can
isotope the entire collection 
$\alpha$ to $\gamma$ using at most $\sum_{i=1}^n a_i = a$ type 2 elementary moves.
Similarly, we can isotop $\beta$ to $\gamma$ using at most $ b $ type 2
elementary moves. It follows that we can isotop $\alpha$ to $\beta$ using
at most $ a+b $ type 2 elementary moves. 

\end{enumerate}

The total number of elementary moves required in these isotopies is at most
twice the number of type 2 elementary moves.

$~~~\qed$

Define the {\em length} of a $PL$ curve transverse to the
1-skeleton of a triangulated surface to be the number of edges of
the 1-skeleton that it crosses, counted with multiplicity.

\begin{lemma} \label{basic.delta}
A basic move in a surface with a triangulation of valence at most $V$ can be realized by at most
28 elementary moves for
a type 1 basic move, $6V+24 $ elementary moves for a type 2 basic move, and 22 elementary moves
for a type 3 basic move. In all cases $17V$ elementary moves suffice. 
\end{lemma}
\bproof
We count the number of elementary moves required for each type of basic move.

Lemma~\ref{polygonal.delta} implies that a type 1 basic move requires at most 14 elementary moves in each of
the two triangles meeting the edge along which the curve is isotoped.

The type 2 basic move can be done in the following sequence of steps. 
\begin{enumerate}
\item Two type 2 elementary moves slide the arc in a neighborhood of an edge to pass through the adjacent
vertex.
\item Two type 2 elementary moves slide the curve
 to run along edges in a neighborhood of the adjacent vertex.
\item $V-2$ type 2 elementary moves slide the arc over a neighborhood of the vertex in each adjacent
triangle.
\item At most 7 type 2 elementary moves slide the arc to basic curves in each of the original two
triangles, by Lemma~\ref{polygonal.delta}. 
\item Two type 2 elementary moves slide the arc to basic curves in each of the $V-2$
triangles where new intersections were introduced. 
\end{enumerate}
Summing, this basic move requires less than $3V + 12$
type 2 elementary moves, and $6V + 24$ in total.

\begin{figure}[hbtp]
\centering
\includegraphics[width=.8\textwidth]{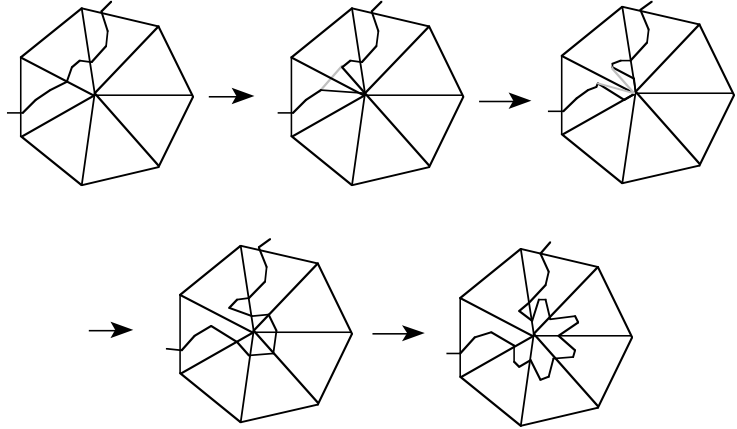}
\caption{ Elementary moves for a type 2 basic move.}
\label{delta.1} 
\end{figure}

The type three basic move requires
\begin{enumerate}
\item One type 2 elementary move to slide the arc to coincide with a subarc of its adjacent edge.
\item At most 10 type 2 elementary moves slide the
resulting arc to a basic curve with the same
endpoints in the adjacent triangle, by an application of Lemma~\ref{polygonal.delta}. This sums to at most 11
elementary moves of type 2 and 22 elementary moves of all types.
\end{enumerate}

Since $V \ge 3$, in all cases $17V$ elementary moves suffice. $~~~\qed$

If two basic curves are isotopic using
only basic moves of type 1, we say that they are {\em parallel} basic curves. 

\begin{lemma} \label{parallel}
Let $\alpha$ and $\beta$ be two parallel basic curves with total length $\ell$.
Then $\alpha$ can be isotoped to coincide with $\beta$ by
$\ell/2$ basic moves of type 1. Similarly, $\alpha$
can be isotoped across
$\beta$ by $\ell/2$ basic moves of type 1. \end{lemma}
\bproof
Each of $\alpha$ and
$\beta$ has length $\ell/2$ and they intersect each triangle of the surface in
pairs of parallel basic arcs.
Using $\ell/2$ basic moves of type 1, we can either move the intersections of 
$\alpha$ along each edge to coincide with the intersection of
$\beta$ with that edge, or move the intersections past to the other side
of $\beta$.
$~~~\qed$

Let $\alpha$ and $\beta$ be intersecting simple closed curves in general position on a surface $F$. A
{\em 2-gon} is a subdisk of $F$ whose boundary consists of one subarc of each of $\alpha$ and
$\beta$. A 2-gon is {\em innermost} if no arc of $\alpha$ or $\beta$ meets its interior. In some
cases, it can be shown that curves with excess intersection create innermost 2-gons. This is
analyzed in some generality in \cite{Hass-Scott85}. The following is Lemma 3.1 of
\cite{Hass-Scott85}.

\begin{lemma} \label{innermost}
Let $F$ be a compact orientable surface and let $\alpha$ and $\beta$ be isotopic simple closed curves
on $F$, intersecting transversely.
Then either $\alpha$ and $\beta$ are disjoint or there is an innermost 2-gon on $F$ bounded by one
arc from each of $\alpha$ and
$\beta$.
\end{lemma}
\bproof
The lift of $\alpha$ to the universal cover of $F$ is a family of disjoint lines $\tilde \alpha$, and the lift of $\beta$ is a similar family $\tilde \beta$. Assume $\alpha$ and $\beta$ intersect. Then a lift $\alpha_1$ of $\alpha$ and a
lift $\beta_1$ of $\beta$ also intersect. The curves $\alpha_1$ and $\beta_1$ are isotopic to
disjoint lines by a lift of an isotopy making $\alpha$ and $\beta$ disjoint. The fundamental group of $F$ acts on the universal cover, and $\alpha_1$ and $\beta_1$ have a
common cyclic stabilizer. It follows that they intersect infinitely often. An arc of each
between adjacent crossing points cuts off a 2-gon in the universal cover. A finite number of arcs
of $\tilde \alpha \cup \tilde \beta$ meet this 2-gon. Each such arc is disjoint from either $\tilde \alpha$ or $\tilde \beta$, so that it cuts off a smaller 2-gon inside the previous one. We can pass to an innermost 2-gon whose interior meets no arc of $\tilde
\alpha \cup \tilde \beta$. This 2-gon is disjoint from any of its translates, and must project 1-1
down to $F$.
$~~~\qed$

\begin{lemma} \label{isotopy.length}
Let $F$ be an orientable triangulated surface with $u$ triangles, and let $\alpha$ and $\beta$ be isotopic simple basic
curves on $F$ intersecting transversely.
Let $l$ be the length of $\alpha \cup \beta$, and let $V$ be the maximal valence
of a vertex in $F$. There is an isotopy carrying $\beta$ to $\alpha$ via no more than $l^4uV$ basic moves.
\end{lemma}
\bproof
We treat cases. Note that $V \ge 3, l \ge 4$ and $u \ge 4$.

{\em Case 1: $\alpha$ and $\beta$ are disjoint.} 

Then $\alpha$ and $\beta$
bound a subsurface of $F$ homeomorphic to an annulus. We will perform a series of basic moves which isotop the curves
to coincide.

The edges of $F$ induces a cellulation of the annulus.
The cells are subdisks of triangles,
bounded by arcs in $\alpha \cup \beta$.
If there is a vertex in this cellulation, we can slide
$\alpha$ or $\beta$ over a vertex by a basic move of type 2.
Repeating at most $u$ times, we eliminate
all vertices in the annulus.
The length of $\alpha \cup \beta$ is increased by at most $V-2$
each time we do a basic move, so the final
curves have total length at most $l+u(V-2) = l+uV -2u$. 

The annulus between the curves is cut by the edges
of $F$ into a collection of cells.
Since $\alpha$ and $\beta$ are disjoint, and each cell boundary contains an edge lying on the 1-skeleton between each edge lying on $\alpha \cup \beta$, all cells are even
sided. An Euler characteristic calculation shows that either all cells are 4-gons, or there
is a 2-gon cell in the annulus, bounded by a basic arc of type 2 and an arc on 
an edge of the triangulation.

If there is a cell which is a 2-gon, we can perform a basic move of type 3 which slides an arc of $\alpha$ or $\beta$
across that 2-gon. This decreases the length of $\alpha \cup \beta$ by two.
After repeating at most $(l+uV - 2u)/2$ times, we
eliminate all 2-gon cells from the annulus, and arrive at a cellulation in which all cells are 4-gons. When all cells are
4-gons, the two boundary components of the annulus are parallel basic curves, each of length at most
$(l+uV-2u)/2$. Using at most $(l+uV-2u)/2$ basic moves of type 1, we can make the two curves
coincide.

The total number of basic moves required in the three steps above is bounded by 
$$
u + (l+uV-2u)/2 + (l+uV-2u)/2 = l +uV. 
$$
{\em Case 2: $\alpha$ and $\beta$ intersect.} 

Two basic arcs in a triangle with
distinct endpoints intersect at most once, so the total
number of intersections of $\alpha$ and $\beta$ is at most $l^2/4$.
Lemma~\ref{innermost} gives the existence of an
innermost disk $D_1$ on the surface, bounded by
subarcs $\alpha_1$ of $\alpha$ and $\beta_1$ of $\beta$.
The 1-skeleton of $F$ induces a cellulation
of $D_1$.
The cells are bounded by arcs in $\alpha_1 \cup \beta_1$ and arcs in the 1-skeleton of $F$. If
$D_1$ contains a vertex of $F$, we can isotop $\alpha_1$ or $\beta_1$ over a vertex (possibly a different one) by a basic
move of type 2, keeping
$\alpha \cup \beta$ embedded. The number of vertices is bounded above by $u$,
so repeating at most $u$ times, we eliminate all vertices in $D_1$, using at most
$u$ basic moves of type 2. The length of
$\alpha_1 \cup \beta_1$ is increased by at most $V-2$ by each basic move, so $\alpha$ and
$\beta $ are isotoped to curves with total length less than $ l+u(V-2) = l + uV - 2u$. 


Since $\alpha_1$ and $\beta_1$ meet at exactly two points in the boundary of $D_1$,
there are at most two cells in $D_1$ with an odd number of sides. An Euler characteristic calculation shows that one of
three cases occurs: 

{\em Case~1.}
All of $D_1$ is contained in a single triangle, or 

{\em Case~2.}
There is a 2-gon cell in $D_1$, or

{\em Case~3.}
All cells of $D_1$ are 4-gons, other then two 3-gons containing points of $\alpha_1 \cap \beta_1$.

The first case is impossible because basic arcs don't intersect in more than one point in a single triangle. A 2-gon cell
in $D_1$ gives a basic arc of type 2. We can use this to perform a basic move of type 3, which slides an arc of
$\alpha_1$ or $\beta_1$ across a 2-gon. This decreases the length of $\alpha_1 \cup \beta_1$ by two, and thus can be
repeated at most $ (l+uV -2u)/2 $ times. This leaves us in Case 3. We can then isotop $\beta_1$ across $\alpha_1$, using
at most $(l+uV-2u)/2$ basic moves of type 1. We have then eliminated the innermost 2-gon $D_1$, reducing the number of
intersection points of $\alpha$ and $\beta$ by two, while using at most $ l + uV $ basic moves.

The length of $\alpha \cup \beta$ is
now bounded by the larger value $l+uV$. If $\alpha$ and $\beta$ still intersect, we repeat with another innermost disk.
Eliminating this disk requires at most $ (l+uV) + uV = l + 2uV$ basic moves. Repeating at most $l^2/4$ times we
eliminate all intersections and arrive at disjoint curves.
Summing the resulting series, which has at most $l^2/4$ terms,
gives a bound of 
$$ l+uV + (l + 2uV ) \dots + (l + (l^2/4)uV ) \le l^3/4 + \frac{l^2}{4}(\frac{l^2}{4}-1) uV/2 \le
l^3/4 + l^4uV/32 $$ 
basic moves. With the curves now disjoint, and having length less than 
$\ds l+\frac{l^2}{4}uV$, case 1 isotops them to coincide using at most $\ds l+ \frac{l^2}{4}uV +uV $ basic moves.
The total number of basic moves used is less than $$l^3/4 + l^4uV/32 + l+ \frac{l^2}{4}uV +uV .$$ 
Since $l \ge 4$ and $u \ge 4$,
$$
l^4uV \ge l^4uV/32 + 2l^3uV \ge l^4uV/32 + l^3/4 + 4l^2uV \ge l^4uV/32 + l^3/4 + l + l^2uV/4 + uV
$$
so the sum is
less than $l^4uV$ as claimed.
$~~~\qed$

\begin{lemma} \label{surface.bound}
Let $F$ be an orientable triangulated surface with $u$ triangles, and let $\alpha$ and $\beta$ be isotopic basic curves on $F$ in
general position. Let $l$ be the total length of $\alpha \cup \beta$. Then there is an isotopy carrying $\beta$ to $\alpha$ via no
more than $ 17 l^4 u^3$ elementary moves. \end{lemma}
\bproof
The required number of basic moves is at most $l^4uV$, and each basic move requires at most $17V$
elementary moves, so the total number of elementary moves required is less than
$17 l^4 u V^2$. We know $V \le u$
so that $ 17 l^4 u^3 $ gives a bound on the number of elementary moves. $~~~\qed$

We now give the proof of the main result of this section. 

\paragraph{Proof of Theorem ~\ref{surface.bound.thm}:} We begin by isotoping $\alpha$ to a basic curve. First we perturb
$\alpha$ to miss the zero-skeleton of the triangulation. There are at most $u$ vertices, and each has at most $V$
edges coming into it, where $V$ is the maximal valence of the triangulation of the surface. We can perturb $\alpha$ near
any vertex which it passes through so that it crosses at most $V/2$ edges near that vertex. 
This requires a total of at most $u(V/2+1)$ type 2 elementary moves, and $u V+2$ elementary moves of both types. The number of
segments on the curve has now been increased to at most $s + uV/2 +u < s+uV$. 

Next we isotop the intersection of $\alpha$ with each triangle to a union of basic arcs. This can increase the number
of segments by at most a factor of three, to $3(s+uV)$ and requires less than $ 4(s+UV)$ elementary
moves by Lemma~\ref{polygonal.delta}. We can do similar moves on $\beta$, resulting in new curves $\alpha$ and
$\beta$ which are transversely intersecting basic curves with a total of at most $ 6(s+uV)$ segments.
Note that the length of a basic curve is between one third and one half of the number of segments it contains,
therefore $\alpha \cup \beta$
has length at most $3(s+uV)$. We now apply Lemma~\ref{surface.bound},
isotoping the curves together using at most $ 17 [3(s + uV)]^4 u^3 $ elementary moves.
We can bound $V$ by $u$,
obtaining a bound of at most
$$ u^2 + 2 + 4(s + u^2) + 1377[ (s + u^2)]^4 u^3 < 2^{11} (s + u^2)^4 u^3$$ 
elementary moves. $ ~~~\qed $

\section{Proof of Theorem~\ref{thm102}}
\hsp
We combine the results of \S2--\S5 to prove Theorem~\ref{thm102}. 

\paragraph{Proof of Theorem~\ref{thm102}:}
To prove Theorem~\ref{thm102}, we first construct a compact triangulated $PL$ submanifold $M_K$ of $M$ obtained by removing a
solid torus neighborhood of $K$. We barycentrically subdivide $M$ twice and construct a new triangulated manifold 
$ M_K = M - R_K ~ ,$ by removing the interior of a regular neighborhood $R_{K}$ of $K$. 
The boundary of $M_K$ consists of a subdivision of the original boundary $\partial M$ plus a 2-torus $\partial R_K$,
the peripheral torus. The manifolds $M_K$ and $R_K$ contain between them $576t$ tetrahedra. As outlined in
\S3, we deform the original knot $K$ across an annulus by elementary moves to a longitude $K_1$ on the peripheral
2-torus $\partial R_K$. Then we deform $K_1$ on the peripheral torus by elementary moves to the boundary $K_2$ of
the embedded normal disk given by Haken's Theorem~3.1.
Finally we deform the curve $K_2$ across this disk using
elementary moves to a triangle.

We count elementary moves.
By Lemma~\ref{disk.count}
the normal disk with $K_2$ as boundary contains at most 
$ 2^{8 \cdot 576t+6} = 2^{4608t+6} $ triangles,
and the number of edges $|K_2|$ of $K_2$ is at most
\begin{equation} 
|K_2| \le 2^{4608t+7} < 2 ^{4609t} = s_2 ~. 
\label{eqN62}
\end{equation}
Also by Lemma~\ref{disk.count} it takes at most $2^{4608t+7}$
elementary moves to transform $K_2$ to a triangle. 

By Theorem~\ref{longitude.count} the longitude $K_1$ consists of at most 
\begin{equation}
|K_1| \le 2^{1858(576t)+1} < 2^{1070209t} = s_1
\label{eqN63} 
\end{equation}
edges, and the isotopy from $K$ to $K_1$
requires at most the same number of elementary moves.

It remains to bound the number of elementary moves needed to carry $K_1$ to $K_2$, using
Theorem~\ref{surface.bound.thm}. We note that $K_1$ is isotopic to $K_2$ if
properly oriented.
Theorem~\ref{surface.bound.thm} applies with $s = \max (s_1, s_2) = s_1$
to yield a bound of at most $ 2^{11} (s + u^2)^4 u^3$
elementary moves taking $K_1$ to $K_2$, where $u$ is the number
of triangles on the peripheral torus.
Since $u \le 576t < 2^{10}t$,
Corollary~\ref{surface.bound.cor} yields a bound of
$ 2^{10^7t-1}$ for the number of elementary moves taking $K_1$ to $K_2$.
For the three steps together we obtain a bound of
$ 2 ^{4609t} + 2^{1070209t} + 2^{10^7t-1} < 2^{10^7t}$.
~~~$\qed$ 

\paragraph{Remark.}
The assumption $K \subseteq Interior (M)$ was used to allow the peripheral torus $R_K$ to be constructible by two
barycentric subdivisions and to have $\partial R_K$ disjoint from $\partial M$. For the general case that $K
\subseteq M$, we could proceed by barycentrically subdividing $M$ twice, then isotoping $K$ off the boundary $\partial M$
to a ``parallel'' knot $K''$ contained strictly in the interior of $M''$, using $O(t)$ elementary moves. Now Theorem
1.2 applies to $K''$ in $(M'', \partial M'' )$.

\section{Reidemeister moves and elementary moves in $\RR^3$} \hsp
The remainder of the paper is concerned with 3-manifolds embedded in $\RR^3$, and is devoted to proving the
Reidemeister move bound of Theorem~\ref{thm101}.

We outline the proof of Theorem~\ref{thm101}.
We have given a knot diagram $\sD$ with $n$ vertices. We use it to construct a triangulated $PL$ manifold $M$
embedded in $\RR^3$ which contains a knot $K$ in its 1-skeleton whose projection on the $xy$-plane in $\RR^3$ is the
knot diagram $\sD$. In section 8 we show that such a manifold $M$ can be constructed using $O(n)$ tetrahedra. Now
Theorem~\ref{thm102}
shows that we can isotop $K$ to a trivial knot using $2^{c_3 n}$ elementary moves. We project these
elementary moves down onto the $xy$-plane and obtain a sequence of Reidemeister moves. The proof is completed by
bounding the number of Reidemeister moves in terms of elementary moves, which we do below. 

In the rest of this section we describe Reidemeister moves and relate them to elementary moves. We recall basic
facts about knot and link diagrams. 

A projection of a link in $\RR^3$ into a plane is {\em regular} if all crossings are in general position, and at most
2~edges intersect the pre-image of any point. Corresponding to a regular projection of a link is a {\em link
diagram}. This is a labeled planar graph having the following properties: 
\begin{enumerate}
\renewcommand{\labelenumi}{(\roman{enumi})} 
\item
Any connected component of $\sD$ with no vertices is an (isolated) loop. 
\item
Each vertex has exactly four adjacent edges, with two edges labeled ``undercrossing'' and two labeled
``overcrossing'' at the vertex.
\item The labeling prescribes a cyclic ordering of edges at each vertex, and ``overcrossing edges'' are adjacent to
``undercrossing edges'' in this ordering.
\end{enumerate}
\noindent
Conversely, for any labeled planar graph $\sD$ satisfying (i)--(iii) there exists an unoriented link $L$ having
projection $\sD$. We note that the cyclic ordering of vertices in (iii) actually specifies a planar embedding of the
graph $\sD$, since it determines the faces of the planar embedding. 

A {\em link component} in a link diagram $\sD$ is obtained by tracing a closed path of edges such that each pair of consecutive edges match at their common vertex, either both overcrossing or both undercrossing. The {\em order} of $\sD$ is the number of components it contains. A {\em link diagram component} of $\sD$ is a
connected component of the underlying graph; it can be the union of several link components.
The {\em crossing measure} $| \sD |$ of a link diagram is 
\begin{equation}
| \sD | = \# (\mbox{vertices}) + \# (\mbox{link diagram components}) - 1 ~. 
\label{eq201}
\end{equation}

A {\em knot diagram} is a link diagram that consists of one link component. A {\em trivial knot diagram}
is a knot diagram that has no vertices. If $\sD$ is a knot diagram, then $| \sD|$ is the {\em
crossing number} of the knot projection. 
\begin{lemma} \label{regular.proj} 
If a polygonal link $L$ has regular projection $\sD$, then
\begin{equation}\label{eq202}
| {\sD} | \leq |L|^2 ~.
\end{equation}
\end{lemma}
\bproof
The number of vertices $|L|$ is at least as large as the number of distinct line segments it contains.
Since line segments project to line segments, each pair of line segments produces at most one crossing.
There are less than $|L|^2$ pairs. Furthermore line segments that project into different connected components of $| {\sD} |$
cannot produce any crossings.
Thus the number of nonintersecting pairs of line segments is at least
$\#(\mbox{components of~} {\sD}) -1$, so
\ref{eq202} follows.
$~~~\qed$

The bound \ref{eq202} is the correct order of magnitude, because one can easily construct
a sequence of knots $K$ with $|K| \ra \infty$, which have the property that $| {\sD} | > \frac{1}{2} |K|^2$.
There is no upper bound for $|L|$ in terms of $| {\sD} |$, but there does exist a link $L$ with a
projection having diagram ${\sD}$, with
$|L| \leq c | {\sD} |$,
see \cite[Lemma~5.1]{HasLagPip}
and \S 3.

\begin{figure}[hbtp]
\centering
\includegraphics[width=.5\textwidth]{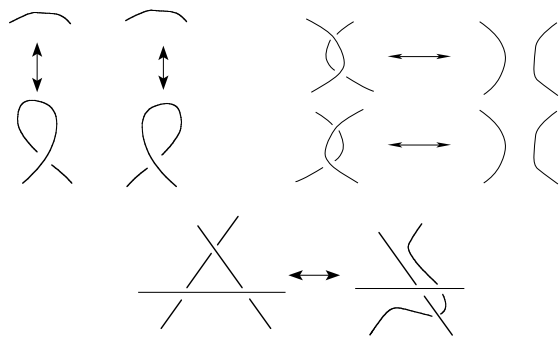}
\caption{Reidemeister moves of types 1, 2 and 3.}
\label{reid.moves}
\end{figure}

Reidemeister moves are local moves on link diagrams that are analogous to elementary moves.
They convert a link diagram ${\sD}$ to another link diagram ${\sD}'$ by a local change in its structure.
For unoriented links there are four kinds of Reidemeister moves, examples of which are pictured in
Figure~\ref{reid.moves}, see one of 
\cite[p.~48]{Mur96}, \cite{Adams},\cite{BurZie85},\cite{},\cite{Rolfsen}.
(For oriented links there are more Reidemeister moves,
see Wu \cite{Wu92}.) 

We call two knot diagrams ${\sD}_1$ and ${\sD}_2$ {\em equivalent} if and only if one can be converted into the other by a finite sequence of Reidemeister moves. This relation corresponds to knot equivalence: Knots $K_1$ and $K_2$ have knot diagrams ${\sD}_1$ and ${\sD}_2$ that are equivalent if and only if ${\sD}_1$ and ${\sD}_2$ are equivalent.

Reidemeister moves on a link diagram ${\sD}$ are related to elementary moves on a link $L$ whose projection is ${\sD}$, but the correspondence of moves is not one-to-one. 

\begin{lemma} \label{reid.2}
Let $L$ be a polygonal link and suppose that it has diagram ${\sD}$ as regular projection on a plane.
Let $L'$ be a link derived from L by an elementary move of type (2) or $(2')$, and suppose
that it has regular projection ${\sD}'$ on the same plane. Then there is a sequence of at most
$2 |L| + 2 | {\sD} |$
Reidemeister moves that converts ${\sD}$ to ${\sD}'$. \end{lemma} 

\bproof
As a line segment is moved affinely in $\RR^3$ with one end fixed, its projection intersects the projection of any vertex of the polygonal link $L$ in at most one point. Therefore it can generate at most one type~I or type~II
Reidemeister move involving the projection of that vertex. The number of vertices is bounded by $|L| $.
Similarly, each pair of crossing segments in ${\sD}$ can produce at most one type~III Reidemeister move.
Since two
line segments are moved in a type (2) or $(2')$ elementary move, this yields at most $2 |L| + 2 | \sD |$
Reidemeister moves. $~~~\qed$ 

The bound of Lemma~\ref{reid.2} is nearly best possible. 

We use these lemmas to prove a bound for Reidemeister moves in terms of elementary moves.

\begin{theorem} \label{Elementary.to.Reidemeister} 
Let $L$ and $L'$ be polygonal links each consisting of at most $n$ line
segments and suppose that ${\sD}$ and
${\sD}'$ are link diagrams obtained by regular projection of $L$ and $L'$ on a fixed plane.
If $L$ can be transformed to $L'$ using $k$ elementary moves, then ${\sD}$ can be transformed to ${\sD}'$ using at most $2k(n+ \frac{1}{2} k +1)^2$ Reidemeister moves. \end{theorem}

{\bf Proof.} Let $L = L_0 , L_1 , L_2 , \ldots , L_k = L'$
be the set of intermediate links occurring in the sequence
of elementary moves, and let ${\sD}_i$ denote the projection of $L_i$ on the fixed plane. We may assume that all
projections are regular, by applying a
a transversality argument and infinitesimally moving the plane of projection \cite[Prop. 1.12]{BurZie85}. Since an elementary move changes the number of vertices by at most
one, and
$|L_0 | , |L_k | \leq n$, we have, for
$1 \leq i \leq k$, that
\begin{equation}
|L_i | \leq \min (n+i , n+k -i ) \leq n + \frac{1}{2} k ~. 
\label{eq203} 
\end{equation} 
Lemma~\ref{regular.proj} gives
\begin{equation} | {\sD}_i | \leq |L_i |^2 \leq (n+ \df{1}{2} k)^2 = n^2 + nk + \df{k^2}{4} ~. \label{eq204}
\end{equation}
Thus we obtain
$$
2| L_i | + 2 | {\sD}_i | \leq 2 (n+ \df{1}{2} k +1 )^2 ~, $$ and applying
Lemma~\ref{reid.2}
to each $L_i$ separately yields the upper bound $2k(n+ \frac{1}{2} k+1)^2$.
$~~~\qed$

\section{Triangulation of knot complements in $\RR^3$} \hsp
In this section we construct a triangulated 3-manifold $M_K$ which is a knot complement in a ball
embedded in $\RR^3$. That is, $M_K = S^3 - (R_K \cup B)$, in which $R_K$ is an open
regular neighborhood of the knot $K$, $B$ is an open ball containing the point at infinity in its
interior, and
$\bar{R}_K \cap \bar{B} = \emptyset$. The manifold $M_K$ has boundary $\partial M_K = \partial R_K \cup
\partial B$ which consists of two connected components, which are topologically a 2-sphere and a
torus.
As in \S2, we
call such a manifold $M_K$ a {\em truncated knot complement}. 

We carry out this construction more generally for link diagrams, as it requires no extra work. Given a link diagram ${\sL}$ with crossing measure $n$, we construct a triangulated manifold $M$ in $\RR^3$ which contains $O(n)$ tetrahedra and has a link $L$ in its 1-skeleton which has a regular projection that gives the link diagram ${\sL}$. 

\begin{theorem} \label{triangulate}
Given a link diagram ${\sL}$ of crossing measure $n$,
one can construct a triangulated convex polyhedron $P$ in $\RR^3$ such that:
\begin{enumerate}
\renewcommand{\labelenumi}{(\roman{enumi})} 
\item The triangulation contains at most $840n$ tetrahedra.
\item Every vertex in the triangulation is a lattice point $(x,y,z) \in \ZZ^3$, with $0 \leq x \leq 30n$,
$0 \leq y < 30n$ and $-6 \leq z \leq 6$.
\item There is a link ${\sL}$ embedded in the
1-skeleton of the triangulation which lies entirely in the interior of $P$, and whose orthogonal
projection on the plane $z = 0$ is regular and is a link diagram isomorphic to ${\sL}$.
\end{enumerate}
\end{theorem}
\bproof
The link diagram ${\sL}$ comes with a topological planar embedding which determines the faces.
We begin by adding extra vertices to
${\sL}$. To each edge or non-isolated loop we add two new vertices of degree~2, which splits it into three edges. To each isolated loop we add three new vertices of degree~2, making it an isolated triangle. The resulting
labeled graph ${\sG}$ is still planar, has no loops or multiple edges, and has at most $5n$ vertices. (The
worst case consists of several disjoint single crossing projections.) Let $m$ denote the number of
vertices of ${\sG}$, and call the $m-n$ vertices added {\em special vertices}. Note that each vertex of the original graph has associated to it four distinct special vertices. 

\begin{figure}[hbtp]
\centering
\includegraphics[width=.4\textwidth]{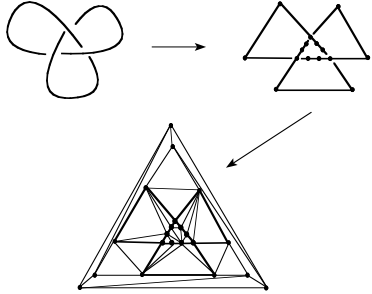}
\caption{Triangulating a knot projection.}
\label{triangulation}
\end{figure}

Using the topological planar embedding of ${\sG}$, we next add extra edges to triangulate each face of ${\sG}$,
obtaining a topologically embedded triangulated planar graph ${\sG}'$. The graph ${\sG}'$ has $m$ vertices and $2m-4$
bounded triangular faces, and the unbounded face is also a triangle. We next encase the graph ${\sG}'$ in a
triangulated planar graph ${\sH}'$ with $m+3$ vertices and $2m+2$ bounded triangular faces, by framing the outside
triangular face with a larger triangle and then adding 6~edges to subdivide the resulting concentric region into
6~triangles, obtaining ${\sG}''$. See Figure~\ref{triangulation}.

Our next goal is to replace the topological embedding of ${\sG}''$ in the plane with a straight-line embedding. In 1952 Fary showed that
straight-line embeddings always exist, and a recent result of de Frajsseix, Pach and Pollack
\cite{deFPacPol90} states that the vertices can be chosen to lie in a small rectangular grid of integer lattice points. Their results imply that there exists a planar embedding of ${\sH}'$ whose vertices
$V({\sH}' )$ lie in the plane $z = -1$, and are contained in the grid 
\begin{displaymath}
\{(x,y,-1): 0 \leq x,y
\leq 10n-1 ~; \quad x,y \in \ZZ \}~, 
\end{displaymath}
and all edges of the graph ${\sH}'$ are straight line segments. 

Next we make an identical copy ${\sH}''$ of ${\sH}'$ on the plane $z=1$ with ${\sH}'' = {\sH}'+ (0,0,2)$, and
vertex set $V({\sH}'') = V({\sH}') + (0,0,2)$. We now consider the polyhedron $P'$ which is the convex
hull of $V({\sH}')$ and $V({\sH}'')$. It is a triangular prism, because the outside face of ${\sH}'$ is
a triangle. We add $m_1$ vertical edges connecting each vertex $\bv \in V({\sH}' )$ to its copy
$\bv + (0,0,1) \in V({\sH}'')$.
Let ${\sE}$ denote these edges, together with all the edges in ${\sH}'$ and ${\sH}''$. Using these edges, the polyhedron $P'$ decomposes into $2m+2$ triangular prisms $\{Q_j : 1 \leq
j \leq 2m-4\}$, with the top and bottom faces of each being congruent triangular faces of ${\sH}'$ and
${\sH}''$.

We next triangulate $P'$ by dissecting each of the $2m+2$ triangular prisms into 14
tetrahedra, as follows. We subdivide each vertical rectangular face into four triangles using its
diagonals. Then we cone each rectangular face to the centroid of $P'$, and note that the centroid lies in
the plane $z = 0$. We add 4 new vertices to each prism, one on each
rectangular face and one in the center. The point
of this subdivision is that the triangulations of
adjacent prisms are compatible. Note also that all new
vertices added lie in the lattice $\frac{1}{3} \ZZ^2$. 

Let ${\sE}'$ denote all the edges in the union of these triangulated prisms. We identify the link diagram ${\sG}$ with the graph ${\sG}$ embedded in ${\sH}'$. We next observe that there is a
polygonal link $L$ embedded in the edge set ${\sE}'$ whose projection is the link diagram
${\sL}$. We insist that any edge in ${\sD}$ that runs to an undercrossing have its vertex
which is an undercrossing be in the plane $z = -1$, while any edge that runs to an
overcrossing have the overcrossing vertex in the plane $z = 1$. Each such edge goes from one of the
$n$ original vertices of ${\sL}$ to a special vertex adjacent to that vertex. Edges that do not meet vertices labeled
overcrossing or undercrossing are assigned to the $z = -1$ plane. The edge corresponding to an edge running from the
$z = -1$ to the $z = 1$ plane is contained in one of the diagonals added to a prism. The resulting embedding $L$ in
$P'$ has a regular orthogonal projection onto the
$(x,y)$-plane, since no edges are vertical and only transverse double points occur as
singularities in the projection.
Note that no edge of $L$ touches the outer triangle of ${\sH}'$ because all edges project to
the graph ${\sG}'$ contained strictly inside ${\sH}'$. This will be important for property (iii).

The knot edges lie in the upper and lower boundary of $P'$, but do not touch its sides. To
have them not touch the boundary we take two additional copies of $P'$ and glue one to its top, along the plane $z =1$, and one to its bottom, along the plane
$z =-1$. The knot now lies entirely in the interior of the resulting polyhedron $P''$. All
the vertices added in this construction lie in $\frac{1}{3} \ZZ^2$. The total number of
tetrahedra used in triangulating $P''$ is $84(m+1)$, which is at most $420(n+1)$, or at most
$840n$ for $n \geq 1$.

Now triple all coordinates $(x,y,z) \rightarrow (3x, 3y, 3z)$ to obtain a triangulated
polytope $P = 3P''$ with integer vertices, with a triangulation that satisfies (i)--(iii).
$~~~\qed$ \\

\paragraph{Remark.}
(1). The construction in Theorem~\ref{triangulate}
can be effectively computed in time $O(n \log n )$
using the Hopcroft-Tarjan planarity algorithms
\cite{HopTar74} and the embedding
algorithm of de Frajsseix, Pach and Pollack \cite[Section~4]{deFPacPol90}.

(2).
One can ask whether a polygonal knot $K$ in $\RR^3$ that consists of $n$ line segments can be
embedded in the 1-skeleton of some triangulated polyhedron that uses $O(n)$ tetrahedra. We do
not know if this can always be done. There exist polygonal knots whose planar projections in any
direction have at least $\Omega (n^2 )$ crossings. Avis and El Gindy \cite{AviElG87}
showed that a set of $n$ points in general position in $\RR^3$ can be embedded as vertices of a
triangulation of their convex hull using $O(n)$ tetrahedra, but that in the general case $\Omega
(n^2 )$ tetrahedra are sometimes needed. 

\section{Proof of Theorem~\ref{thm101}}
\hsp
We combine the earlier results to prove Theorem~\ref{thm101} with an explicit bound on the number of
Reidemeister moves required to isotop an unknot to the trivial projection. 

\paragraph{Proof of Theorem~\ref{thm101}.}
Applying Theorem~\ref{triangulate}
to the knot diagram ${\sD}$ yields a triangulated
convex polyhedron $M$ embedded in $\RR^3$ with
$t$ tetrahedra, where,
\begin{displaymath} 
t \le 840n ~,
\end{displaymath}
and with a knot $K$ embedded in the 1-skeleton of the interior of $M$ that projects to ${\sD}$.
The hypotheses of Theorem~\ref{thm102}
are satisfied for $M$, hence there is a sequence of at
most $2^{c_2 t}$ elementary moves
transforming $K$ to the trivial knot in $M$, with $c_2 = 10^7$.
Now Theorem~\ref{Elementary.to.Reidemeister} 
gives a bound of at most
\begin{equation}
2^{c_2 t} \left(6 \cdot 840n + \frac{1}{2} \cdot 2^{c_2 t} +1
\right)^2 ~.
\label{eqN92} 
\end{equation} 
Reidemeister moves, since $K$ contains at most $6 \cdot 840n$ line segments. 

The bound (\ref{eqN92}) yields the desired bound of the form $2^{c_1 n}$ for an explicit constant
$c_1$. Using the value $c_2 = 10^7$
from the proof of Theorem~\ref{thm102}, one finds that $c_1 =
10^{11}$ suffices.
~~~$\qed$

\vspace*{.5in}
\begin{tabular}{ll}
{\tt email}: & {\tt hass@math.ucdavis.edu} \\ & {\tt jcl@research.att.com} \\
~~~ \\
~~~ \\
{\tt address}: & Prof. J. Hass \\
& Dept. of Mathematics \\
& Univ. of California, Davis \\
& Davis, CA 95616 \\
~~~ \\
& Dr. J. C. Lagarias \\
& Room C235, AT\&T Labs--Research \\
& 180 Park Avenue \\
& Florham Park, NJ 07932-0971
\end{tabular}

\end{document}